\documentclass{cocv}
\usepackage{fix-cm}
\usepackage[scr=rsfs]{mathalpha}
\usepackage{amsmath,amssymb,amsfonts}
\usepackage{graphicx}
\DeclareMathOperator{\tr}{tr}

\begin{document}

\title{Finite-Horizon LQR for General Markov Jump Linear Systems: Deterministic Reformulation and Reduced-Order Solution}
\thanks{The authors thank the Vice-Rectory for Research, Extension and Social Projection (VIEPS) of Universidad del Atlántico for support.}

\author{Alfredo R. R. Narváez}\address{Department of Mathematics, Universidad del Atlántico, Barranquilla, Colombia; e-mail: alfredoroa@uniatlantico.edu.co}
\author{Jeinny Peralta}\address{Department of Mathematics, Universidad del Atlántico, Barranquilla, Colombia; e-mail: jeinnyperalta@uniatlantico.edu.co}
\author{M. A. C.  Candezano}\address{Department of Mathematics, Universidad del Atlántico, Barranquilla, Colombia; e-mail: miguelcaro@uniatlantico.edu.co}

\date{Received ...; Revised ...; Accepted ...}

\begin{abstract}
This paper studies the Linear Quadratic Regulator (LQR) problem for continuous-time Markov Jump Linear Systems (MJLS) governed by general finite-state Markov chains that may include transient, absorbing, or non-communicating states. The problem, posed over a finite time horizon, is reformulated deterministically by expressing the cost functional in terms of a collection of second-moment matrices of the system state. A projection operator is introduced to restrict the analysis to the subspace of visited states, namely those with positive probability of being reached within the time horizon. The solution of the resulting deterministic problem is obtained from a reduced-order system of coupled matrix Riccati differential equations involving only the visited states, which define a quadratic value function satisfying a Hamilton–Jacobi–Bellman-type equation. The structure of this equation is formally justified in the matrix setting via the Riesz–Fréchet representation theorem, establishing a rigorous foundation for the deterministic reformulation and resolving an open aspect in previous literature. Several numerical examples, including cases with non-communicating states, validate the theoretical results and illustrate the practical relevance of the proposed generalization.
\end{abstract}



\maketitle

\section{Introduction}

It is common to observe that many natural phenomena undergo abrupt changes in their behavioral dynamics; that is, phenomena whose \emph{modus operandi} changes at each instant and for which it is not possible to have prior knowledge of when or how these changes will occur, even though there may be some information about their behavior in each mode of operation. A phenomenon or process that exhibits this type of behavior is known in the literature as a stochastic process. Of course, it is also possible to fit such a process into a mathematical model describing another phenomenon, for example, via a system of differential equations, such that the dynamics of the model parameters are governed by the behavior of the stochastic process. According to the literature, the equations that describe such models are known as jump dynamical systems; and in particular, when the system evolves linearly in its state variable and the parameter jumps are governed by a Markov chain, the system is referred to as a Markov Jump Linear System (MJLS), which can be written as:
\begin{equation}\label{mjls}
\begin{aligned}
    \dot{x}(t) &= A_{\theta(t)}(t)x(t) + B_{\theta(t)}(t)u(t), \\
    x(0) &= x_0,
\end{aligned}
\end{equation}
where, for $t \geq 0$, $x(t) \in \mathbb{R}^n$ denotes the state trajectory of the system with initial condition $x_0$; $u(t) \in \mathbb{R}^m$ is a piecewise continuous function referred to as the control law applied to the system; the parameters $A_{\theta(t)}(t)$ and $B_{\theta(t)}(t)$ are piecewise continuous functions belonging to finite families of matrix-valued functions $\{A_1(t),\ldots,A_N(t)\}$ and $\{B_1(t),\ldots,B_N(t)\}$, respectively; and $\{\theta(t)\}_{t\geq 0}$ is a homogeneous continuous-time Markov chain with finite state space $\mathcal{S} = \{1,\ldots,N\}$ and generator matrix $\Lambda \in \mathbb{R}^{N \times N}$ with entries $\lambda_{i,j}$ for $i,j = 1,\ldots,N$. Such a process is defined on a probability space $(\Omega, \mathcal{F}, \text{Pr})$ equipped with a filtration $\{\mathcal{F}_t\}_{t \geq 0}$ such that $\theta(t)$ is $\mathcal{F}_t$-measurable for all $t \geq 0$. Here, $\Omega$ denotes the sample space, $\mathcal{F}$ is a $\sigma$-algebra, $\text{Pr}$ is a probability measure, and $\{\mathcal{F}_t\}_{t \geq 0}$ is a collection of sub-$\sigma$-algebras of $\mathcal{F}$ satisfying $\mathcal{F}_s \subseteq \mathcal{F}_t$ for all $s < t$. The process $\{\theta(t)\}_{t\geq 0}$ satisfies the Markov property in the sense that 
\[
\text{Pr}(\theta(t_{n+1})=j \mid \theta(t_n)=i, \theta(t_{n-1})=i_{n-1}, \ldots, \theta(t_0)=i_0) = \text{Pr}(\theta(t_{n+1})=j \mid \theta(t_n)=i),
\]
for all $j, i, i_{n-1}, \ldots, i_0 \in \mathcal{S}$ and $t_{n+1} > t_n > \cdots > t_0 \geq 0$. The distribution of the Markov chain is characterized by $\phi(i) := \text{Pr}(\theta(0) = i)$, $i \in \mathcal{S}$, and $P_{i,j}(t) := \text{Pr}(\theta(s + t) = j \mid \theta(s) = i)$, $i, j \in \mathcal{S}$, which represent the initial state distribution and the transition probabilities, respectively. The matrix $P(t) = (P_{i,j}(t))$ for $i,j \in \mathcal{S}$ and $t \geq 0$ is called the transition probability matrix, and it is stochastic in the sense that $\sum_{j \in \mathcal{S}} P_{i,j}(t) = 1$ for all $t \geq 0$ and $i \in \mathcal{S}$. The transition rates $\lambda_{i,j}$ can be obtained from $P(t)$ via the limits $\lim_{t \rightarrow 0^+} \frac{P_{i,i}(t) - 1}{t} = \lambda_{i,i}$ and $\lim_{t \rightarrow 0^+} \frac{P_{i,j}(t)}{t} = \lambda_{i,j}$ for $i \neq j$. These rates satisfy $\lambda_{i,j} \geq 0$ for $i \neq j$, and $\lambda_{i,i} = -\sum_{j \neq i} \lambda_{i,j}$. In addition, the initial condition $x_0$ is generally assumed to be a random variable independent of the process $\{\theta(t)\}_{t \geq 0}$, with mean $\bar{x}_0$ and covariance matrix $\Sigma$.

The problem addressed in this paper lies at the intersection of stochastic systems and applied optimal control. The formulation and analysis developed here are motivated not only by theoretical interest but also by the wide range of applications in engineering, economics, and systems subject to regime-switching dynamics. In this sense, the methods employed contribute to the broader field of applied mathematics and align with ongoing research efforts that seek to bridge stochastic modeling, control theory, and functional analysis. Given its relevance in numerous applications, optimal control theory has naturally extended to the framework of MJLS. In this context, the present work focuses on studying certain optimal control problems for such systems, specifically the so-called finite-horizon LQR problem with jumps 
defined as:
\begin{equation}\label{lqr:intro}
    \begin{aligned}
    \min_{u(t) \in \mathcal{U}} & \;\; \mathscr{E}\left[\int_{0}^{T}\left(x(t)'Q_{\theta(t)}(t)x(t) + u(t)'R_{\theta(t)}(t)u(t)\right) dt + x(T)'Q_{\theta(T)}(T)x(T)\right] \\
    & \hspace{-.4cm} \text{s.t.} \quad \begin{cases}
    \dot{x}(t) = A_{\theta(t)}(t)x(t) + B_{\theta(t)}(t)u(t), \\
    x(0) = x_0
    \end{cases}
    \end{aligned}
\end{equation}
that is, to find a control input $u(t)$, for $0 \leq t \leq T$, within a set $\mathcal{U} \subseteq \mathbb{R}^m$ of admissible control inputs, that minimizes the cost functional in~\eqref{lqr:intro}, where $\mathscr{E}(\cdot)$ denotes the expectation operator; $(\cdot)'$ is the matrix transpose; and $Q_{\theta(t)}(t)$ and $R_{\theta(t)}(t)$ are given semi-positive and positive definite matrices, respectively, for all $0 \leq t \leq T$, with $T > 0$ fixed. 

Unlike the vast majority of the existing literature on optimal control for MJLS, including the seminal works \cite{Wonham73,Mariton,fragoso1991optimal,Ji92,fragoso2010separation,ValleCosta2013}, we do \emph{not} assume that the Markov chain $\{\theta(t)\}_{t\geq 0}$ is irreducible or that all states communicate. That is, we allow for the general case where the chain may contain transient states, absorbing states, or multiple non-communicating classes. Formally, we do not require that for any pair of distinct states $i,j \in \mathcal{S}$, there exist $s,t > 0$ for which $P_{i,j}(s) > 0$ and $P_{j,i}(t) > 0$. This generalization is motivated by practical applications where certain operational modes may only be accessible from specific initial conditions, or where failure modes act as absorbing states. 

On the other hand, in recent years, the study of optimal control for MJLS has gained significant attention due to its applicability to stochastic dynamic systems with multiple operational modes. Contemporary works have addressed the LQR problem under various assumptions, including model-free approaches \cite{Fan2025}, distributed control in power systems \cite{Qin2025}, and infinite-horizon cases with ergodic risk constraints \cite{TalebiLi2025}. Predictive control and policy optimization methods have further advanced the field \cite{SooJean2023,JanschPorto2020}. However, these studies, while addressing diverse problem settings, also depend on Markov chains satisfying irreducibility or ergodicity assumptions, thus excluding systems with transient, absorbing, or multiple non-communicating classes of states. The rigorous extension of LQR theory to arbitrary finite-state Markov chains, particularly the formal derivation of the HJB equation and systematic treatment of non-visited states, remains an open problem. This work contributes to filling that gap by providing a deterministic reformulation and reduction framework that naturally accommodates these general conditions. These recent advances motivate the present deterministic reformulation, which accommodates general Markov chain structures and provides a rigorous Hilbert-space framework for the finite-horizon LQR problem.


We begin by reviewing the solution to the finite-horizon LQR problem with jumps using typical tools from control theory such as dynamic programming, while accounting for the stochastic nature of the problem, as presented in, for example, \cite{Wonham73,Mariton,fragoso1991optimal,Ji92,fragoso2010separation,ValleCosta2013}. These classical results, however, fundamentally rely on the assumption that all states of the Markov chain communicate, that is, the chain is irreducible. Under this restrictive assumption, conditional expectations of the form $\mathscr{E}(x(t)x(t)' \mid \theta(t) = i)$ are well-defined for all states $i \in \mathcal{S}$ and all times $t \in [0,T]$, which is essential for formulating the coupled Riccati differential equations that characterize the optimal control. However, since solving a deterministic problem is generally less complex than addressing its stochastic counterpart, and motivated by the need to relax the communicating states assumption, the purpose of this work is to reformulate the finite-horizon LQR problem with jumps using deterministic optimal control results, suitably adapted to a new state space of a different nature than that of the original system state.

We introduce as the new state variable a collection $\mathbf{X}(t)$ of $N$ second-moment matrices of the trajectory $x(t)$, whose $i$-th component is defined as $X_i(t) = \mathscr{E}\left(x(t)x(t)' \cdot 1_{\{\theta(t) = i\}}\right)$ for $i \in \mathcal{S}$, where $1_{\{\cdot\}}$ denotes the indicator function. Crucially, this definition remains well-posed even when the Markov chain is not irreducible. For states $i$ that are never visited, i.e., states with $\text{Pr}(\theta(t) = i) = 0$ for all $t \in [0,T]$, we have $X_i(t) = 0$ by definition, thus avoiding the need for conditional expectations on zero-probability events. The collection $\mathbf{X}(t)$ takes values in the Hilbert space $\mathbb{H}_{\mathbb{R}}^{n+}$ consisting of tuples $\mathbf{V} = (V_1, \ldots, V_N)$ of symmetric, positive semidefinite $n \times n$ real matrices, i.e., $\mathbf{V} \in \mathbb{H}_{\mathbb{R}}^{n+}$ if and only if $V_i = V_i' \geq 0$ for each $i$. This space is a subspace of the larger Hilbert space $\mathbb{H}_{\mathbb{R}}^n$, whose elements are collections $\mathbf{U} = (U_1, \ldots, U_N)$ of real $n \times n$ matrices. 

To properly account for states that may never be visited, we introduce the notion of visited states $\mathscr{Z} \subseteq \mathcal{S}$, defined as those states $i$ for which $\text{Pr}(\theta(t) = i) > 0$ for some $t \in [0,T]$. We then define a projection operator $\mathcal{P}: \mathbb{H}_{\mathbb{R}}^n \to \mathbb{H}_{\mathbb{R}}^n$ that zeros out components corresponding to non-visited states, and introduce the subspace $\mathbb{V} := \text{Im}(\mathcal{P})$ of collections that vanish on non-visited states. Since $\mathbf{X}(t) \in \mathbb{V}$ for all $t \in [0,T]$, the optimization problem naturally reduces to this lower-dimensional subspace. Assuming the minimizing control takes the form $u(t) = L_{\theta(t)}(t)x(t)$, and introducing suitable restricted operators $\mathcal{K}_{\mathbf{U}}^{\mathscr{Z}}: \mathbb{V} \to \mathbb{V}$ for a given parameter $\mathbf{U} \in \mathbb{H}_{\mathbb{R}}^n$ that act only on the visited states, we follow and extend the approach in \cite{Roa&ECosta} to reformulate the finite-horizon LQR problem with jumps as:
\begin{equation}\label{lqr_X}
    \begin{aligned}
    \min_{\mathbf{L}} & \int_0^T \langle \mathbf{\Gamma}(t, \mathbf{L}), \mathbf{X}(t) \rangle dt + \langle \mathbf{\Gamma}(T), \mathbf{X}(T) \rangle \\
    \text{s.t.} \quad & \begin{cases}
    \dot{\mathbf{X}}(t) = \mathcal{K}_{\mathbf{A+BL}}^{\mathscr{Z}}(\mathbf{X}(t)), \\
    \mathbf{X}(0) := \mathbf{X}^0 \in \mathbb{V},
    \end{cases}
    \end{aligned}
\end{equation}
where $X_i(0) = \mathscr{E}(x(0)x(0)' \cdot 1_{\{\theta(0) = i\}}) = \phi(i) x_0 x_0'$ for $i \in \mathscr{Z}$, and $X_i(0) = 0$ for $i \notin \mathscr{Z}$; $\Gamma_i(t,\mathbf{L}) := Q_i(t) + L_i'(t) R_i(t) L_i(t)$ for $0 \leq t < T$, and $\Gamma_i(T) := Q_i(T)$. We then adapt tools from the Dynamic Programming Algorithm to derive the following Hamilton--Jacobi--Bellman (HJB) equation in the Hilbert space setting:
\begin{equation}\label{HJB_meq}
\begin{aligned}
0 = \min_{\mathbf{L}} \Big\{ \langle \mathbf{\Gamma}(t, \mathbf{L}), \mathbf{X} \rangle + \partial_t J(t, \mathbf{X}) + \langle \nabla_{\mathbf{X}} J(t, \mathbf{X}), \mathcal{K}_{\mathbf{A+BL}}^{\mathscr{Z}}(\mathbf{X}) \rangle \Big\}, \\
J(T, \mathbf{X}) = \langle \mathbf{\Gamma}(T), \mathbf{X} \rangle, \quad \forall \mathbf{X} \in \mathbb{V}.
\end{aligned}
\end{equation}

Although the approach developed here follows the general line of reasoning in \cite{Roa&ECosta}, that earlier work did not provide a formal derivation of the HJB equation in terms of second-moment matrix collections. Instead, the structure of the HJB equation was assumed without rigorous justification, and the origin of such an equation in the Hilbert space setting remained unclear. In contrast, the present work makes two main contributions. First, we establish the mathematical foundation for the HJB formulation by invoking the Riesz--Fréchet representation theorem, which provides a rigorous justification for the form of the HJB equation in the context of Hilbert spaces of matrix-valued collections. This closes an important conceptual gap left open in \cite{Roa&ECosta} and enhances the theoretical soundness of the deterministic reformulation. Second, building upon the generalization initiated in \cite{Roa&ECosta}, we provide a comprehensive treatment of the LQR problem for arbitrary finite-state Markov chains by systematically introducing the projection operator $\mathcal{P}$ and restricted operators that naturally handle non-visited states. This allows us to derive a reduced-order system of coupled Riccati equations involving only the $|\mathscr{Z}|$ visited states, rather than all $N$ states, and to fully characterize the solution structure in this general setting.

This conceptual clarification, in turn, enables the formulation of a sufficiency verification theorem that establishes a fundamental connection between the HJB equation in \eqref{HJB_meq} and the minimization of the original cost functional in \eqref{lqr:intro}. Specifically, under suitable regularity conditions, we show that if a continuously differentiable function $V(t, \mathbf{X})$ satisfies the HJB equation in \eqref{HJB_meq} and the associated minimizer $\mathbf{L}^*(t)$ exists for all $t \in [0,T]$, then $V$ coincides with the value function and the feedback law $\mathbf{L}^*(t)$ solves the finite-horizon LQR problem with jumps. Furthermore, we prove that this minimizer $\mathbf{L}^*(t)$ can be computed explicitly by solving a system of $|\mathscr{Z}|$ coupled matrix Riccati differential equations, one for each visited state. In fact, the HJB equation admits a solution of the form $V(t, \mathbf{X}) = \langle \mathbf{Y}(t), \mathbf{X} \rangle$, where $\mathbf{Y}(t) \in \mathbb{V}$ solves the reduced-order Riccati system. Hence, the coupled Riccati equations provide both necessary and sufficient conditions for solving the original stochastic control problem, even in the presence of non-communicating states.

The organization of this work is as follows. Section 2 presents a review of the classical finite-horizon LQR problem with jumps under the irreducibility assumption, introduces the necessary mathematical framework including the projection operator and restricted operators, and develops the deterministic reformulation in terms of second-moment matrices that naturally handles general Markov chains. Section 3 develops the solution of the reformulated problem using deterministic optimal control techniques, including the rigorous derivation of the HJB equation via the Riesz--Fréchet representation theorem and the reduced-order coupled Riccati system. Section 4 provides illustrative numerical examples that validate the theoretical results, including cases with non-communicating states that demonstrate the practical relevance of our generalization. One example involves an application that models the optimal trajectory design required to keep a satellite in orbit despite potential failures in its thrusters, where failure modes are modeled as absorbing states. Finally, Section 5 contains the conclusions derived from the results presented.

\section{Preliminaries}

In this section, we present a review of the finite-horizon LQR problem with Markovian jumps. We begin by recalling its classical formulation in the stochastic setting under the restrictive assumption that all states of the underlying Markov chain communicate with each other (i.e., irreducibility). We then introduce its deterministic reformulation based on the second-moment matrix framework, which serves as the foundation for removing this assumption and developing a general theory applicable to arbitrary finite-state Markov chains.

\subsection{Overview of the LQR Problem with Markovian Jumps}

The formulation of a finite-horizon LQR problem for a MJLS was first introduced in \cite{Wonham73}. This formulation, also found in \cite{Mariton} and later used in \cite{Ji92} to solve a linear-quadratic-Gaussian (LQG) problem for a MJLS with additive noise via the separation principle, consists of solving the optimization problem:
\begin{equation}\label{lqr:cond}
\begin{aligned}
\min_{u(t) \in \mathcal U} & \;\; \mathscr{E}\left[\int_{0}^{T}\left(x(t)'Q_{\theta(t)}(t)x(t)+u(t)'R_{\theta(t)}(t)u(t)\right)dt \,\big|\, x(0)=x_0,\theta(0)=i_0\right] \\
& \hspace{-.4cm} \text{s.t.:} \; \begin{cases}
\dot{x}(t)=A_{\theta(t)}(t)x(t)+B_{\theta(t)}(t)u(t), \\
x(0)=x_{0} \in \mathbb{R}^{n}.
\end{cases}
\end{aligned}
\end{equation}

The approach to solving \eqref{lqr:cond} relies on a stochastic version of the Dynamic Programming Algorithm, from which an HJB-type equation is derived:
\begin{align}\label{eq:HJBM}
0 = \min_{u\in \mathbb{R}^{m}}\left\{ x'Q(t)x+u'R(t)u+ \partial_{t} V(t,x,i)+\mathcal{L}_{u}V(t,x,i) \right\},
\end{align}
for all $(t,x,i) \in [0,T]\times \mathbb R^n \times \mathcal S$, with terminal condition $V(x,i,T)=0$. The operator $\mathcal L_u$, acting on the solution $V$ of the HJB equation \eqref{eq:HJBM}, is the infinitesimal generator of the joint Markov process $\{(x(t),\theta(t))\}_{t\geq 0}$ under the admissible input $u$ and is given by:
\begin{equation}\label{gen:inf}
\mathcal L_uV(x,i,t)=\nabla_xV(t,x,i)'(A_i(t)x+B_i(t)u)+\sum_{j=1}^{N}\lambda_{ij}V(t,x,j).
\end{equation}

The result providing the solution to problem \eqref{lqr:cond} is as follows.

\begin{thrm}[\cite{Mariton}]\label{LQRS:Mariton}
The solution to problem \eqref{lqr:cond} is given by:
$$u^*(t)=-R_{i}^{-1}(t)B_{i}'(t)Y_{i}(t)x(t),$$
whenever $\theta(t)=i$, $t\in[0,T]$, where the functions $Y_i$, $i \in \mathcal S$, are obtained from the system of coupled matrix Riccati differential equations:
\begin{equation}\label{Ric:YT=0}
\begin{aligned}
-\dot{Y}_{i}(t)&=A_{i}'(t)Y_{i}(t)+Y_{i}(t)A_{i}(t)+\sum_{j\in S} \lambda_{ij}Y_{j}(t) +Q_{i}(t)-Y_{i}(t)B_{i}(t)R^{-1}_{i}(t)B_{i}'(t)Y_{i}(t), \\
Y_i(T)&=0, \qquad i\in \mathcal S.
\end{aligned}
\end{equation}
The minimum cost is:
$$ J^*=x_0'Y_{i_0}(0)x_0, $$
given that $x(0)=x_0$ and $\theta(0)=i_0$.
\end{thrm}

It is important to note that the cost functional in \eqref{lqr:cond} is computed with respect to the conditional probability distribution of the joint initial pair $(x(0), \theta(0))$ and, unlike the formulation considered in this work, it omits a terminal cost depending on the state at $T$. In contrast, the quadratic cost functional considered here includes it explicitly as the unconditional expectation:
\begin{align}\label{eq:LQRSM}
\mathscr{E}\left[\int_{0}^{T}\left(x(t)'Q_{\theta(t)}(t)x(t)+u(t)'R_{\theta(t)}(t)u(t)\right)dt + x(T)'Q_{\theta(T)}(T)x(T)\right].
\end{align}

Thus, the corresponding finite-horizon LQR problem with jumps is formulated as:
\begin{equation}\label{lqr_saltos}
\begin{aligned}
\min_{u(t) \in \mathcal U} & \;\; \mathscr{E}\left[\int_{0}^{T}\left(x(t)'Q_{\theta(t)}(t)x(t)+u(t)'R_{\theta(t)}(t)u(t)\right)dt + x(T)'Q_{\theta(T)}(T)x(T)\right] \\
& \hspace{-.4cm} \text{s.t.:} \; \begin{cases}
\dot{x}(t)=A_{\theta(t)}(t)x(t)+B_{\theta(t)}(t)u(t), \\
x(0)=x_{0} \in \mathbb{R}^{n},
\end{cases}
\end{aligned}
\end{equation}
where $Q_{\theta(T)}(T) \geq 0$ is the terminal cost weight matrix. The solution to this problem was also reported in \cite{fragoso2010separation} and can be obtained through the separation principle applied to an LQG problem involving a MJLS with additive noise. There, the optimal control law is given by:
\begin{equation}\label{u_op}
u^*(t)=-R_{\theta(t)}^{-1}(t)B_{\theta(t)}'(t)Y_{\theta(t)}(t)x(t), \quad t\in [0,T],
\end{equation}
where $Y_1(t),\ldots,Y_N(t)$ solve the terminal value problem given by the system of coupled matrix Riccati differential equations:
\begin{equation}\label{Ric_eq}
\begin{aligned}
-\dot{Y}_{i}(t)&=A_{i}'(t)Y_{i}(t)+Y_{i}(t)A_{i}(t)+\sum_{j\in S} \lambda_{ij}Y_{j}(t) +Q_{i}(t)-Y_{i}(t)B_{i}(t)R^{-1}_{i}(t)B_{i}'(t)Y_{i}(t), \\
Y_i(T)&=Q_i(T).
\end{aligned}
\end{equation}

In fact, \cite[Theorem 1]{fragoso2010separation} proves the existence of a unique collection of $N$ positive semidefinite matrices $Y_1(t),\ldots,Y_N(t)$, $t \in [0,T]$, solving \eqref{Ric_eq}. The result is summarized as follows.

\begin{thrm}[Adapted from \cite{fragoso2010separation}]\label{lqr:C&F}
Consider the optimal stochastic control problem \eqref{lqr_saltos}. Then the optimal control is given by \eqref{u_op}, and the optimal cost is
\begin{equation}\label{cost:opt:C&F}
J^* = \sum_{i \in \mathcal{S}} \left( \tr(\phi(i)\Sigma_i Y_i(0)) + \phi(i) \bar{x}_0' Y_i(0) \bar{x}_0 \right),
\end{equation}
where $\Sigma_i := \mathscr{E}((x(0) - \bar{x}_0)(x(0) - \bar{x}_0)' \mid \theta(0) = i)$, and $Y_i(t)$, $i \in \mathcal{S}$, are the unique positive semidefinite solutions to the system \eqref{Ric_eq}.
\end{thrm}

\begin{rmrk}\label{Obs:Jmin}
If we define $X_i(0):=\mathscr E\left(x(0)x(0)'1_{\{\theta(0)=i\}}\right)$, then $X_i(0)=\phi(i)(\Sigma_i+\bar x_0\bar x_0')$ for all $i \in \mathcal S$. Therefore, evaluating \eqref{cost:opt:C&F} gives $J^*=\sum_{i \in \mathcal S}\tr(X_i(0)Y_i(0))$. If $x(0)$ is deterministic and equal to $x_0 \in \mathbb R^n$, then $J^*=\sum_{i \in \mathcal S}\phi(i)x_0'Y_i(0)x_0$.
\end{rmrk}

The proof of Theorem \ref{lqr:C&F}, extracted from \cite[Appendix]{fragoso2010separation}, relies on a sequence of auxiliary results. However, it does not explicitly invoke classical tools of optimal control such as dynamic programming or a verification theorem. After defining $V(t):=\mathscr E(x(t)'Y_{\theta(t)}(t)x(t))$, it is shown that:
\begin{equation*}
    \begin{aligned}
    \mathscr{E}\Bigg[\int_{0}^{T}&\left(x(t)'Q_{\theta(t)}(t)x(t)+u(t)'R_{\theta(t)}(t)u(t)\right)dt+x(T)'Q_{\theta(T)}(T)x(T)\Bigg]\\
    &=V(0)+\int_{0}^{T}\mathscr E\left((u(t)+L_{\theta(t)}(t)x(t))'R_{\theta(t)}(t)(u(t)+L_{\theta(t)}(t)x(t))\right)dt,
    \end{aligned}
\end{equation*}
from which the positive definiteness of $R_{\theta(t)}(t)$ naturally gives the optimal control $u^*(t)=-L_{\theta(t)}(t)x(t)$, with $L_{\theta(t)}(t)=R_{\theta(t)}^{-1}(t)B_{\theta(t)}'(t)Y_{\theta(t)}(t)$. Then, using Remark \ref{Obs:Jmin}, it follows that $V(0) = J^*$. 

All the results presented above, from \cite{Wonham73, Mariton, fragoso2010separation}, implicitly assume that the underlying Markov chain $\{\theta(t)\}_{t\ge0}$ is irreducible, meaning that all states communicate with each other. Under this assumption, for any initial distribution $\phi$ with $\phi(i)>0$ for all $i\in\mathcal S$, it follows that $\pi_i(t)=\Pr(\theta(t)=i)>0$ for all $t>0$ and all $i\in\mathcal S$. This ensures that the conditional expectations $\mathscr E(x(t)x(t)'\mid\theta(t)=i)$ are well defined for all states and times, which is essential for the validity of the coupled Riccati equations~\eqref{Ric_eq}. However, in many practical applications the Markov chain may contain transient or absorbing states, or multiple non-communicating classes. For instance, fault detection systems often model distinct failure modes as absorbing states, and hierarchical control systems may involve modes that are only accessible from specific initial conditions. The main contribution of this work is to extend the LQR theory for MJLS to the general case where the irreducibility assumption is removed, allowing for arbitrary finite-state Markov chains. In the next section, we revisit the finite-horizon LQR problem with jumps from a different perspective that naturally accommodates this generalization.

\subsection{Reformulation of the LQR Problem with Jumps}

To set the stage for this alternative treatment of the problem that will enable us to handle general Markov chains, we introduce a reformulation of the LQR problem in which the main variables are elements of a Hilbert space of finite collections of matrices. We begin by establishing the necessary mathematical framework.

\subsubsection{Hilbert Space of Matrix Collections}

Let us consider the set $\mathbb{H}_{\mathbb{R}}^{n,m}$ of all collections of $N$ matrices $\mathbf{V} = (V_1, \ldots, V_N)$, where each $V_i$ belongs to the vector space of all $m \times n$ matrices, denoted by $\mathbb{B}(\mathbb{R}^{n}, \mathbb{R}^{m})$ (with the notations $\mathbb{H}_{\mathbb{R}}^n := \mathbb{H}_{\mathbb{R}}^{n,n}$ and $\mathbb{B}(\mathbb{R}^{n}, \mathbb{R}^{n}) := \mathbb{B}(\mathbb{R}^{n})$ for brevity). The space $\mathbb{H}_{\mathbb{R}}^{n,m}$ is also a vector space, equipped with the operations:
\[
\mathbf{V} + \mathbf{W} := (V_1 + W_1, \ldots, V_N + W_N), \quad \text{and} \quad \alpha \mathbf{V} := (\alpha V_1, \ldots, \alpha V_N),
\]
for all $\mathbf{V} = (V_1, \ldots, V_N), \mathbf{W} = (W_1, \ldots, W_N) \in \mathbb{H}_{\mathbb{R}}^{n,m}$ and any scalar $\alpha \in \mathbb{R}$. Moreover, for $\mathbf{V}, \mathbf{W} \in \mathbb{H}_{\mathbb{R}}^{n,m}$ as above, we write
\[
\mathbf{V} \mathbf{W} := (V_1 W_1, \ldots, V_N W_N).
\]

Additionally, we define the subspaces $\mathbb{H}_{\mathbb{R}}^{n*}$ and $\mathbb{H}_{\mathbb{R}}^{n+}$ as
\[
\mathbb{H}_{\mathbb{R}}^{n*} := \left\{ \mathbf{V} = (V_1, \ldots, V_N) \in \mathbb{H}_{\mathbb{R}}^n \; : \; V_i = V_i', \, \forall i = 1,\ldots, N \right\}
\]
and
\[
\mathbb{H}_{\mathbb{R}}^{n+} := \left\{ \mathbf{V} = (V_1, \ldots, V_N) \in \mathbb{H}_{\mathbb{R}}^{n*} \; : \; V_i \geq 0, \, \forall i = 1,\ldots, N \right\}.
\]

Moreover, we introduce the following bilinear mapping $\langle \cdot \, ; \cdot \rangle : \mathbb{H}_{\mathbb{R}}^n \times \mathbb{H}_{\mathbb{R}}^n \rightarrow \mathbb{R}$, defined by
\begin{align}\label{prointer}
\langle \mathbf{V}; \mathbf{W} \rangle := \sum_{i=1}^N \tr(V_i' W_i), \quad \forall \: \mathbf{V}, \mathbf{W} \in \mathbb{H}_{\mathbb{R}}^n,
\end{align}
which satisfies the axioms of an inner product on $\mathbb{H}_{\mathbb{R}}^n$.

\begin{rmrk}\label{H:Hilbert}
Note that the space $\mathbb{B}(\mathbb{R}^{n})$ is isomorphic to the Euclidean space $\mathbb{R}^{n^2}$ via the isomorphism $\textbf{{\upshape vec}}:\mathbb{B}(\mathbb{R}^{n}) \rightarrow \mathbb{R}^{n^2}$ defined as
\[
\textbf{{\upshape vec}}(A) = \begin{bmatrix} a_{11} & \cdots & a_{n1} & a_{12} & \cdots & a_{n2} & \cdots & a_{1n} & \cdots & a_{nn} \end{bmatrix}', \quad \forall A \in \mathbb{B}(\mathbb{R}^{n}),
\]
with $A = (a_{ij})$ for $i,j = 1, \ldots, n$, see \cite{ValleCosta2013}. Moreover, it is easy to verify that $\textbf{{\upshape vec}}(A)' \textbf{{\upshape vec}}(B) = \tr(A'B)$ for all $A, B \in \mathbb{B}(\mathbb{R}^{n})$, where $\tr(A'B)$ defines the inner product between the matrices $A$ and $B$ in $\mathbb{B}(\mathbb{R}^{n})$. Hence, equipping $\mathbb{B}(\mathbb{R}^{n})$ with the Frobenius norm $\|A\|_{\text{F}}^2 := \langle A, A \rangle_F := \tr(A'A)$ turns it into a Banach space, or equivalently, $\mathbb{B}(\mathbb{R}^{n})$ becomes a Hilbert space when endowed with the inner product $\langle \cdot, \cdot \rangle_F$. This, in turn, implies that the space $\mathbb{H}_{\mathbb{R}}^{n}$, equipped with the inner product defined in \eqref{prointer}, is also a Hilbert space.
\end{rmrk}

\subsubsection{Visited States and Projection Operator}

Before proceeding with the characterization of the second-moment matrices, we introduce the notion of visited states and a projection operator that will play a crucial role in our general formulation.

\begin{dfntn}\label{def:visited_states}
We define the set of \emph{visited states} over the time horizon $[0,T]$ as:
\[
\mathscr{Z} := \left\{ i \in \mathcal{S} : \pi_i(t) > 0 \text{ for some } t \in [0,T] \right\},
\]
where $\pi_i(t) = \text{Pr}(\theta(t) = i)$. 
\end{dfntn}

States $i \notin \mathscr{Z}$ are never visited in the sense that $\pi_i(t) = 0$ for all $t \in [0,T]$. By standard properties of continuous-time Markov chains, if $\pi_i(t) > 0$ for some $t > 0$, then $\pi_i(s) > 0$ for all $s \geq t$. Therefore, $\mathscr{Z}$ depends only on the initial distribution $\phi$ and the generator $\Lambda$ and it can be characterized as consisting of all states $i$ for which there exists a state $j$ with $\phi(j) > 0$ and a finite sequence $i_0 = j, i_1, \ldots, i_k = i$ such that
$
\lambda_{i_0, i_1} \lambda_{i_1, i_2} \cdots \lambda_{i_{k-1}, i_k} > 0.
$

\begin{rmrk}
When the Markov chain is irreducible and $\phi(i) > 0$ for all $i \in \mathcal{S}$, we have $\mathscr{Z} = \mathcal{S}$. This is precisely the setting covered by the classical results in Theorems \ref{LQRS:Mariton} and \ref{lqr:C&F}. In the general case, $\mathscr{Z}$ may be a proper subset of $\mathcal{S}$, and we write $|\mathscr{Z}|$ to denote the cardinality of the set of visited states.
\end{rmrk}

\begin{dfntn}\label{def:projection}
Let ${\mathcal P}:\mathbb{H}_{\mathbb{R}}^{n}\rightarrow \mathbb{H}_{\mathbb{R}}^{n}$ denote the projection operator defined componentwise by:
\begin{equation}\label{def-mathcalP}
({\mathcal P}(\mathbf{U}))_i= \begin{cases}
U_i, & i \in \mathscr{Z}, \\
0, & i\notin\mathscr{Z},
\end{cases}
\end{equation}
for all $\mathbf{U} \in \mathbb{H}_{\mathbb{R}}^{n}$. We denote the image of $\mathcal P$ by 
\[
\mathbb{V} := \text{Im}(\mathcal P) = \left\{ \mathbf{U} \in \mathbb{H}_{\mathbb{R}}^{n} : U_i = 0 \text{ for all } i \notin \mathscr{Z} \right\} \subset \mathbb{H}_{\mathbb{R}}^{n}.
\]
Note that $\mathcal P$ is indeed an orthogonal projection since:
\begin{enumerate}
\item $\mathcal{P}^2 = \mathcal{P}$ (idempotence),
\item $\mathcal{P}$ is self-adjoint: $\langle \mathcal{P}(\mathbf{U}); \mathbf{V} \rangle = \langle \mathbf{U}; \mathcal{P}(\mathbf{V}) \rangle$ for all $\mathbf{U}, \mathbf{V} \in \mathbb{H}_{\mathbb{R}}^{n}$.
\end{enumerate}
Furthermore, if ${\mathcal I}:\mathbb{H}_{\mathbb{R}}^{n}\rightarrow \mathbb{H}_{\mathbb{R}}^{n}$ stands for the identity operator, any $\mathbf{U} \in \mathbb{H}_{\mathbb{R}}^{n}$ admits the orthogonal decomposition
$$\mathbf{U} = {\mathcal P}(\mathbf{U}) + (\mathcal I-\mathcal P)(\mathbf{U}),$$
where ${\mathcal P}(\mathbf{U}) \in \mathbb{V}$ and $(\mathcal I-\mathcal P)(\mathbf{U}) \in \mathbb{V}^{\perp}$, with
\[
\mathbb{V}^{\perp} = \left\{ \mathbf{U} \in \mathbb{H}_{\mathbb{R}}^{n} : U_i = 0 \text{ for all } i \in \mathscr{Z} \right\}.
\]
\end{dfntn}

\subsubsection{Operators on Matrix Collections}

Now, for a given element $\mathbf{U} \in \mathbb{H}_{\mathbb{R}}^n$, we introduce the following operators acting on $\mathbb{H}_{\mathbb{R}}^n$: $\mathcal{K}_{\mathbf{U}}: \mathbb{H}_{\mathbb{R}}^n \to \mathbb{H}_{\mathbb{R}}^n$ and $\mathcal{H}_{\mathbf{U}}: \mathbb{H}_{\mathbb{R}}^n \to \mathbb{H}_{\mathbb{R}}^n$, defined componentwise by:
\begin{align}\label{op:dinX}
\mathcal{K}_{\mathbf{U}, i}(\mathbf{Q}) := U_i Q_i + Q_i U_i' + \sum_{j \in \mathcal{S}} \lambda_{ji} Q_j, \quad \forall\, \mathbf{Q} \in \mathbb{H}_{\mathbb{R}}^n,\; i \in \mathcal{S},
\end{align}
\begin{align}\label{op:dinXt}
\mathcal{H}_{\mathbf{U}, i}(\mathbf{Q}) := U_i' Q_i + Q_i U_i + \sum_{j \in \mathcal{S}} \lambda_{ij} Q_j, \quad \forall\, \mathbf{Q} \in \mathbb{H}_{\mathbb{R}}^n,\; i \in \mathcal{S}.
\end{align}

\begin{rmrk}
The operators defined in \eqref{op:dinX} and \eqref{op:dinXt} are adjoint to each other with respect to the inner product \eqref{prointer} (see \cite[Lemma 3.5]{ValleCosta2013}). That is, for all $\mathbf{Q}, \mathbf{P} \in \mathbb{H}_{\mathbb{R}}^n$,
\[
\langle \mathcal{K}_{\mathbf{U}}(\mathbf{Q}); \mathbf{P} \rangle = \langle \mathbf{Q}; \mathcal{H}_{\mathbf{U}}(\mathbf{P}) \rangle.
\]
This adjointness is a key property in deriving the solution to the finite-horizon LQR problem with jumps, as will be shown in subsequent sections.
\end{rmrk}

Since we are working with general Markov chains where not all states may be visited, it is natural to introduce restricted versions of these operators that act only on the subspace $\mathbb{V}$ of collections whose components vanish for non-visited states.

\begin{dfntn}\label{def:restricted_operators}
For a given element $\mathbf{U} \in \mathbb{H}_{\mathbb{R}}^n$, we define the \emph{restricted operators} $\mathcal{K}_{\mathbf{U}}^{\mathscr{Z}}: \mathbb{V} \to \mathbb{V}$ and $\mathcal{H}_{\mathbf{U}}^{\mathscr{Z}}: \mathbb{V} \to \mathbb{V}$ by:
\begin{align}\label{op:dinX_restricted}
\mathcal{K}_{\mathbf{U}}^{\mathscr{Z}}(\mathbf{Q}) := \mathcal{P}\left(\mathcal{K}_{\mathbf{U}}(\mathbf{Q})\right), \quad \forall\, \mathbf{Q} \in \mathbb{V},
\end{align}
\begin{align}\label{op:dinXt_restricted}
\mathcal{H}_{\mathbf{U}}^{\mathscr{Z}}(\mathbf{Q}) := \mathcal{P}\left(\mathcal{H}_{\mathbf{U}}(\mathbf{Q})\right), \quad \forall\, \mathbf{Q} \in \mathbb{V}.
\end{align}
Equivalently, in component form, for $i \in \mathscr{Z}$:
\begin{align}
\mathcal{K}_{\mathbf{U}, i}^{\mathscr{Z}}(\mathbf{Q}) &:= U_i Q_i + Q_i U_i' + \sum_{j \in \mathscr{Z}} \lambda_{ji} Q_j, \label{op:dinX_comp}\\
\mathcal{H}_{\mathbf{U}, i}^{\mathscr{Z}}(\mathbf{Q}) &:= U_i' Q_i + Q_i U_i + \sum_{j \in \mathscr{Z}} \lambda_{ij} Q_j, \label{op:dinXt_comp}
\end{align}
and for $i \notin \mathscr{Z}$:
\[
\mathcal{K}_{\mathbf{U}, i}^{\mathscr{Z}}(\mathbf{Q}) = \mathcal{H}_{\mathbf{U}, i}^{\mathscr{Z}}(\mathbf{Q}) = 0.
\]
\end{dfntn}

\begin{rmrk}\label{rem:restricted_ops_properties}
The restricted operators have the following important properties:
\begin{enumerate}
\item \textbf{Invariance of $\mathbb{V}$:} They map $\mathbb{V}$ into itself. That is, if $\mathbf{Q} \in \mathbb{V}$, then $\mathcal{K}_{\mathbf{U}}^{\mathscr{Z}}(\mathbf{Q}), \mathcal{H}_{\mathbf{U}}^{\mathscr{Z}}(\mathbf{Q}) \in \mathbb{V}$. This is immediate from the definition since they are obtained by composing with $\mathcal{P}$.

\item \textbf{Adjointness:} The operators $\mathcal{K}_{\mathbf{U}}^{\mathscr{Z}}$ and $\mathcal{H}_{\mathbf{U}}^{\mathscr{Z}}$ remain adjoint to each other when restricted to $\mathbb{V}$: for all $\mathbf{Q}, \mathbf{P} \in \mathbb{V}$,
\[
\langle \mathcal{K}_{\mathbf{U}}^{\mathscr{Z}}(\mathbf{Q}); \mathbf{P} \rangle = \langle \mathbf{Q}; \mathcal{H}_{\mathbf{U}}^{\mathscr{Z}}(\mathbf{P}) \rangle.
\]
This follows from the fact that $\mathcal{P}$ is self-adjoint and the adjointness of $\mathcal{K}_{\mathbf{U}}$ and $\mathcal{H}_{\mathbf{U}}$.

\item \textbf{Reduced sums:} In component form, the sums in \eqref{op:dinX_comp}--\eqref{op:dinXt_comp} run only over $j \in \mathscr{Z}$ rather than over all $j \in \mathcal{S}$. This reflects the fact that for $\mathbf{Q} \in \mathbb{V}$, we have $Q_j = 0$ for all $j \notin \mathscr{Z}$, so the terms with $j \notin \mathscr{Z}$ contribute nothing to the sums. This observation is crucial: it shows that the dynamics of visited states decouple from non-visited states.

\item \textbf{Reduction to classical case:} When the Markov chain is irreducible and all states have positive initial probability (i.e., $\mathscr{Z} = \mathcal{S}$), we have $\mathcal{P} = \mathcal{I}$ and $\mathbb{V} = \mathbb{H}_{\mathbb{R}}^{n}$, so the restricted operators coincide with the original ones:
\[
\mathcal{K}_{\mathbf{U}}^{\mathscr{Z}} = \mathcal{K}_{\mathbf{U}}, \quad \mathcal{H}_{\mathbf{U}}^{\mathscr{Z}} = \mathcal{H}_{\mathbf{U}}.
\]
Thus, our general framework recovers the classical setting as a special case.

\item \textbf{Computational efficiency:} In practical implementations, when $|\mathscr{Z}| \ll N$, working with the restricted operators allows us to solve a reduced-order system involving only $|\mathscr{Z}|$ matrix equations rather than $N$ equations. This can lead to significant computational savings.
\end{enumerate}
\end{rmrk}

\subsubsection{Second-Moment Matrices}

Let us define a collection $\mathbf{X}(t) = (X_1(t), \ldots, X_N(t))$ of $N$ second-moment matrices of the state trajectory $x(t)$, where each component is given at time $t \geq 0$ by:
\begin{align}\label{eq_secmom}
X_i(t) = \mathscr{E}(x(t) x(t)' \cdot 1_{\{\theta(t) = i\}}), \quad i \in \mathcal{S}.
\end{align}

\begin{rmrk}\label{rem:conditional_interpretation}
The definition \eqref{eq_secmom} is well-posed for all $i \in \mathcal{S}$ and all $t \in [0,T]$. When $\pi_i(t) > 0$ (i.e., $i \in \mathscr{Z}$), the matrix $X_i(t)$ can be equivalently written as
$$X_i(t) = \pi_i(t) \mathscr{E}(x(t)x(t)' \mid \theta(t) = i),$$
relating it to the conditional second moment. This conditional expectation is well-defined precisely because $\pi_i(t) > 0$. However, when $\pi_i(t) = 0$ (i.e., $i \notin \mathscr{Z}$), the conditional expectation $\mathscr{E}(x(t)x(t)' \mid \theta(t) = i)$ is not uniquely defined, as we are conditioning on an event of probability zero. Nevertheless, the definition \eqref{eq_secmom} remains well-posed and yields $X_i(t) = 0$ for all $t \in [0,T]$, since the event $\{\theta(t) = i\}$ has probability zero. 
This observation has an important consequence: $\mathbf{X}(t) \in \mathbb{V}$ for all $t \in [0,T]$. That is, the collection of second-moment matrices naturally lives in the subspace $\mathbb{V}$ of collections that vanish on non-visited states. Equivalently, $\mathbf{X}(t) = \mathcal{P}(\mathbf{X}(t))$ for all $t \in [0,T]$.
\end{rmrk}

When the system in \eqref{mjls} is a homogeneous MJLS, i.e., when the dynamics reduce to
\begin{equation}\label{hom_mjls}
\dot{x}(t) = A_{\theta(t)}(t)x(t),
\end{equation}
the collection of second-moment matrices defined in \eqref{eq_secmom} evolves according to the following result.

\begin{thrm}[\cite{CosFra}, adapted]\label{thm:second_moment_dynamics}
For a homogeneous MJLS as in \eqref{hom_mjls}, the collection $\mathbf{X}(t) \in \mathbb{V}$ satisfies, for all $t > 0$,
\begin{align}\label{ecdif}
\dot{\mathbf{X}}(t) = \mathcal{K}_{\mathbf{A}(t)}^{\mathscr{Z}}(\mathbf{X}(t)).
\end{align}
Equivalently, in component form, for all $i \in \mathscr{Z}$,
\begin{align}\label{ecdif_comp}
\dot{X}_i(t) = A_i(t) X_i(t) + X_i(t) A_i'(t) + \sum_{j \in \mathscr{Z}} \lambda_{ji} X_j(t),
\end{align}
while for $i \notin \mathscr{Z}$, we have $X_i(t) = 0$ for all $t \in [0,T]$.
\end{thrm}

\begin{proof}
For $i \in \mathscr{Z}$, the result follows from \cite[Proposition 4.4(b)]{CosFra}, which establishes that
\[
\dot{X}_i(t) = A_i(t) X_i(t) + X_i(t) A_i'(t) + \sum_{j \in \mathcal{S}} \lambda_{ji} X_j(t).
\]
However, since $X_j(t) = 0$ for all $j \notin \mathscr{Z}$ (by Remark \ref{rem:conditional_interpretation}), the sum reduces to $\sum_{j \in \mathscr{Z}} \lambda_{ji} X_j(t)$, yielding \eqref{ecdif_comp}.

For $i \notin \mathscr{Z}$, since $\pi_i(t) = 0$ for all $t \in [0,T]$, the definition \eqref{eq_secmom} immediately yields $X_i(t) = 0$, and consequently $\dot{X}_i(t) = 0$.

Therefore, $\mathbf{X}(t) \in \mathbb{V}$ for all $t \in [0,T]$, and the compact form \eqref{ecdif} emphasizes that the dynamics evolve entirely within the subspace $\mathbb{V}$ of visited states. The use of the restricted operator $\mathcal{K}_{\mathbf{A}(t)}^{\mathscr{Z}}$ ensures that the right-hand side also belongs to $\mathbb{V}$, making the equation consistent.
\end{proof}

\begin{rmrk}
The form \eqref{ecdif_comp} reveals an important structural property: the evolution of $X_i(t)$ for $i \in \mathscr{Z}$ depends only on $X_j(t)$ for $j \in \mathscr{Z}$. In other words, the second-moment dynamics for visited states decouple completely from non-visited states. This justifies working with a reduced-order system of dimension $|\mathscr{Z}|$ rather than the full system of dimension $N$.
\end{rmrk}

\subsubsection{Deterministic Reformulation of the LQR Problem}

Next, we rewrite the cost functional in \eqref{eq:LQRSM} in terms of the second-moment matrix $\mathbf X(t)$, assuming a priori that the minimizing control input is of the form $u(t) = L_{\theta(t)}(t)x(t)$ for $t \in [0,T]$. In fact, we claim that the finite-horizon LQR problem with jumps in \eqref{lqr_saltos} can be equivalently rewritten as:
\begin{equation}\label{lqr_X2}
\begin{aligned}
\min_{\mathbf L} & \int_{0}^{T} \langle \mathbf{\Gamma}(t,\mathbf L);\mathbf{X}(t) \rangle dt + \langle \mathbf{\Gamma}(T);\mathbf{X}(T) \rangle\\
& \hspace{-0.5cm} \text{s.t.:} \; \begin{cases}
\dot{\mathbf{X}}(t) = \mathcal{K}_{\mathbf{A}(t)+\mathbf B(t)\mathbf L(t)}^{\mathscr{Z}}(\mathbf{X}(t)), \\
\mathbf{X}(0):= \mathbf X^0 \in \mathbb{V},
\end{cases}
\end{aligned}
\end{equation}
where the components of the initial collection $\mathbf X^0 \in \mathbb{V}$ are given by 
\[
X_i(0) = \mathscr{E}(x(0)x(0)'1_{\{\theta(0)=i\}}) = \phi(i)x_0x_0', \quad i \in \mathscr{Z},
\]
and $X_i(0) = 0$ for $i \notin \mathscr{Z}$. The variable $\mathbf{L}(t) = (L_1(t), \ldots, L_N(t)) \in \mathbb{H}_{\mathbb{R}}^{m,n}$, defined for $t \in [0,T]$, is the decision variable of the problem, and the collection $\mathbf{\Gamma}(t,\mathbf{L}) \in \mathbb{H}_{\mathbb{R}}^{n}$ is defined componentwise as 
\[
\Gamma_i(t,\mathbf{L}) := Q_i(t) + L_i'(t)R_i(t)L_i(t), \quad 0 \leq t < T, \quad i \in \mathcal{S},
\]
and $\Gamma_i(T) := Q_i(T)$ for each $i \in \mathcal{S}$.

We now verify that the cost functional in \eqref{eq:LQRSM} coincides with that in \eqref{lqr_X2}. 
For the terminal term, using the definition \( X_i(T) = \mathscr{E}(x(T)x(T)'1_{\{\theta(T)=i\}}) \), we have
\begin{align}\label{terminal_cost_derivation_compact}
\mathscr{E}\!\left[x(T)'Q_{\theta(T)}(T)x(T)\right] 
& =\sum_{i\in \mathcal{S}} \tr \left( \mathscr{E}\left[x(T)x(T)'1_{\{ \theta(T)=i\}} \right]Q_i(T) \right) \\
& =\sum_{i\in \mathcal{S}} \tr (X_i(T)Q_i(T)) = \langle \mathbf{\Gamma}(T);\mathbf{Q}(T) \rangle,
\end{align}
and since \(\mathbf{\Gamma}(T)=\mathbf{Q}(T)\), this equals \(\langle \mathbf{\Gamma}(T),\mathbf{X}(T)\rangle\). Similarly, if \(u(t)=L_{\theta(t)}(t)x(t)\) for \(t\in[0,T)\), then 
\[
x(t)'(Q_{\theta(t)}(t)+L_{\theta(t)}'(t)R_{\theta(t)}(t)L_{\theta(t)}(t))x(t)
= \sum_{i\in\mathcal S} x(t)' \Gamma_i(t,\mathbf L)x(t)\,1_{\{\theta(t)=i\}},
\]
and hence, by linearity of expectation and the definition of \(X_i(t)\),
\begin{align}\label{integral_cost_derivation_compact}
\mathscr{E}\!\left[\int_0^T x(t)'(Q_{\theta(t)}(t)+L_{\theta(t)}'(t)R_{\theta(t)}(t)L_{\theta(t)}(t))x(t)\,dt\right]
&= \int_0^T\!\sum_{i\in\mathcal S}\tr(\Gamma_i(t,\mathbf L)X_i(t))\,dt \notag\\
&= \int_0^T\!\langle\mathbf{\Gamma}(t,\mathbf L),\mathbf{X}(t)\rangle\,dt.
\end{align}
Thus, the stochastic cost functional \eqref{eq:LQRSM} admits the deterministic representation \eqref{lqr_X2}. 
Since \(X_i(t)=0\) for all \(i\notin\mathscr Z\) and \(t\in[0,T]\) (see Remark~\ref{rem:conditional_interpretation}), the inner products effectively reduce to the visited states:
\[
\langle\mathbf{\Gamma}(t,\mathbf L),\mathbf{X}(t)\rangle
= \sum_{i\in\mathscr Z}\tr(\Gamma_i(t,\mathbf L)X_i(t)).
\]
This observation is crucial: it shows that the optimal control problem naturally decouples into a reduced-order problem over the visited states $\mathscr{Z}$, which may be strictly smaller than the full state space $\mathcal{S}$ when the Markov chain is not irreducible. Moreover, since the constraint $\dot{\mathbf{X}}(t) = \mathcal{K}_{\mathbf{A}(t)+\mathbf B(t)\mathbf L(t)}^{\mathscr{Z}}(\mathbf{X}(t))$ preserves the property that $\mathbf{X}(t) \in \mathbb{V}$ (as the restricted operator maps $\mathbb{V}$ into itself), the entire optimization problem can be formulated over the reduced-order space $\mathbb{V}$ of dimension $n^2 |\mathscr{Z}|$ rather than the full space $\mathbb{H}_{\mathbb{R}}^{n}$ of dimension $n^2 N$. To formalize the consistency between the full-order and reduced-order formulations when control is present, we establish the following result.

\begin{prpstn}[Consistency of restricted and full systems]\label{prop:system_consistency}
Consider the full-order system of controlled dynamics
\begin{equation}\label{eq:full_system}
\dot{\mathbf{X}}(t) = \mathcal{K}_{\mathbf{A}(t)+\mathbf{B}(t)\mathbf{L}(t)}(\mathbf{X}(t)), \quad \mathbf{X}(0) = \mathbf{X}_0 \in \mathbb{H}_{\mathbb{R}}^{n},
\end{equation}
and the restricted system
\begin{equation}\label{eq:restricted_system}
\dot{\mathbf{X}}(t) = \mathcal{K}_{\mathbf{A}(t)+\mathbf{B}(t)\mathbf{L}(t)}^{\mathscr{Z}}(\mathbf{X}(t)), \quad \mathbf{X}(0) = \mathcal{P}(\mathbf{X}_0).
\end{equation}
Then, if $\mathbf{X}_{\text{full}}(t)$ solves \eqref{eq:full_system} and $\mathbf{X}_{\text{res}}(t)$ solves \eqref{eq:restricted_system}, we have
\[
\mathcal{P}(\mathbf{X}_{\text{full}}(t)) = \mathbf{X}_{\text{res}}(t), \quad \forall t \in [0,T].
\]
In other words, the projection of the solution to the full system coincides with the solution to the restricted system.
\end{prpstn}

\begin{proof}
Define $\widetilde{\mathbf{X}}(t) := \mathcal{P}(\mathbf{X}_{\text{full}}(t))$. Clearly, $\widetilde{\mathbf{X}}(0) = \mathcal{P}(\mathbf{X}_0)$. Since $\mathcal{P}$ is a linear and time-independent operator,
\[
\frac{d}{dt}\widetilde{\mathbf{X}}(t) = \mathcal{P}\left(\dot{\mathbf{X}}_{\text{full}}(t)\right) = \mathcal{P}\left(\mathcal{K}_{\mathbf{A}(t)+\mathbf{B}(t)\mathbf{L}(t)}(\mathbf{X}_{\text{full}}(t))\right).
\]
By Definition~\ref{def:restricted_operators}, for any $\mathbf{Q}\in\mathbb{H}_{\mathbb{R}}^{n}$,
\[
\mathcal{K}^{\mathscr Z}_{\mathbf{A}(t)+\mathbf{B}(t)\mathbf{L}(t)}(\mathbf{Q}) = \mathcal{P}\big(\mathcal{K}_{\mathbf{A}(t)+\mathbf{B}(t)\mathbf{L}(t)}(\mathbf{Q})\big).
\]
Writing $\mathbf{Q}=\mathcal{P}(\mathbf{Q})+(\mathcal{I}-\mathcal{P})(\mathbf{Q})$ with $\mathcal{P}(\mathbf{Q})\in\mathbb{V}$ and $(\mathcal{I}-\mathcal{P})(\mathbf{Q})\in\mathbb{V}^\perp$, and noting that for $i\in\mathscr Z$ we have $\lambda_{ji}=0$ whenever $j\notin\mathscr Z$, it follows that
\[
\mathcal{P}\!\left(\mathcal{K}_{\mathbf{U}}(\mathbf{Q})\right)
= \mathcal{P}\!\left(\mathcal{K}_{\mathbf{U}}(\mathcal{P}(\mathbf{Q}))\right)
= \mathcal{K}_{\mathbf{U}}^{\mathscr Z}(\mathcal{P}(\mathbf{Q})).
\]
Applying this to $\mathbf{Q}=\mathbf{X}_{\text{full}}(t)$ gives
\[
\frac{d}{dt}\widetilde{\mathbf{X}}(t)
= \mathcal{K}^{\mathscr Z}_{\mathbf{A}(t)+\mathbf{B}(t)\mathbf{L}(t)}(\widetilde{\mathbf{X}}(t)),
\]
which shows that $\widetilde{\mathbf{X}}(t)$ satisfies \eqref{eq:restricted_system}. By uniqueness of solutions to ODEs, $\widetilde{\mathbf{X}}(t)=\mathbf{X}_{\text{res}}(t)$ for all $t\in[0,T]$.
\end{proof}

\begin{rmrk}\label{rem:control_gains_unvisited}
An important question arises: what about the control gains $L_i(t)$ for $i \notin \mathscr{Z}$? Since states $i \notin \mathscr{Z}$ are never visited, the values of $L_i(t)$ for such states do not affect the cost functional or the state dynamics. They can be chosen arbitrarily (or set to zero for definiteness). The optimization problem \eqref{lqr_X2} effectively determines only $L_i(t)$ for $i \in \mathscr{Z}$.
\end{rmrk}

The reformulated version, expressed in terms of the matrix-valued trajectory $\mathbf{X}(t) \in \mathbb{V}$ and the operator-based structure with explicit use of the restricted operators, provides a rigorous foundation for extending the LQR theory to general Markov chains without the communicating states assumption. This framework will serve as the basis for deriving the main results in the following sections, where we will obtain a reduced-order system of coupled Riccati differential equations involving only the visited states $\mathscr{Z}$.

\section{Solution to the Finite-Horizon LQR Problem with Jumps: A Deterministic Approach}\label{Cap4}

In this section, we present the deterministic approach to solving the finite-horizon LQR problem with Markovian jumps, previously reformulated in terms of second-moment matrices. A key advantage of this reformulation is that it naturally accommodates general Markov chains without the restrictive assumption of irreducibility. By working within the subspace $\mathbb{V}$ of collections that vanish on non-visited states and employing the restricted operators introduced in Section 2, we derive a reduced-order system of coupled Riccati differential equations involving only the visited states $\mathscr{Z} \subseteq \mathcal{S}$. By discretizing the time interval and applying dynamic programming principles, we derive a HJB-type equation associated with the cost-to-go functional. Under suitable regularity assumptions, we first establish the differentiability of the value function and validate the sufficiency theorem that ensures that a solution to the HJB equation also solves the original optimization problem.

\subsection{A Sufficient Condition via the HJB Equation}

We consider a discretization of the problem in \eqref{lqr_X2}, to which, under certain conditions, we apply the dynamic programming algorithm in order to obtain an HJB-type equation. To this end, let us take a partition $\{0=t_0, t_1,\ldots, t_\ell=T\}$ of the interval $[0,T]$ and assume $t_{k+1}-t_k=\delta>0$ for all $k=0,\ldots,\ell-1$, so that $t_k=k\delta$, for $k=0,\ldots, \ell$. We also denote $\mathbf X(t_k)=\mathbf X(k\delta)$ by $\mathbf X^k$ and $\mathbf L(t_k)=\mathbf L(k\delta)$ by $\mathbf L^k$. With this, the dynamics of the second-moment matrices is approximated by:
\[
\mathbf X^{k+1}=\mathbf X^k+\delta\mathcal{K}_{\mathbf{A}^k+\mathbf B^k\mathbf L^k }^{\mathscr{Z}}(\mathbf{X}^k),
\]
where we recall that $\mathbf{X}^k \in \mathbb{V}$ (i.e., components corresponding to non-visited states remain zero), and the cost functional to be minimized, by finding suitable collections of matrices $\mathbf L^0,\ldots,\mathbf L^{\ell-1}$ of appropriate dimensions, is approximated by:
\begin{equation}\label{cost:disc}
\sum_{k=0}^{\ell-1}\langle \mathbf{\Gamma}(t_{k},\mathbf {L}^{k});\mathbf{X}^{k} \rangle \delta + \langle \mathbf{\Gamma}(t_\ell);\mathbf{X}^\ell \rangle.
\end{equation}
Now, for any time $t$ in $[0,T]$ and any collection $\mathbf X \in \mathbb{V}$, we define the cost-to-go associated with problem in \eqref{lqr_X2} as:
\begin{equation}\label{cost:cont}
    J^*(t,\mathbf X)=\min_{\mathbf L(\tau),\, \tau \in [t,T] }\left\{\int_t^T \langle \mathbf{\Gamma}(\tau,\mathbf {L}(\tau));\mathbf{X}(\tau) \rangle d\tau + \langle \mathbf{\Gamma}(T);\mathbf{X}(T) \rangle\right\},
\end{equation}
where the minimization is understood over $\mathbf{X}(\tau) \in \mathbb{V}$ satisfying the constraint $\dot{\mathbf{X}}(\tau) = \mathcal{K}_{\mathbf{A}(\tau)+\mathbf{B}(\tau)\mathbf{L}(\tau)}^{\mathscr{Z}}(\mathbf{X}(\tau))$ with $\mathbf{X}(t) = \mathbf{X}$, and whose discretized version, for any $t_k$ and $\mathbf X^k=\mathbf X$, is given by:
\begin{equation*}
    \tilde J^*(t_k,\mathbf X)=\min_{\mathbf L^k,\ldots,\mathbf L^{\ell-1}}\left\{\sum^{\ell-1}_{i=k} \langle \mathbf{\Gamma}(t_{i},\mathbf {L}^{i});\mathbf{X}^{i} \rangle \delta + \langle \mathbf{\Gamma}(t_\ell);\mathbf{X}^\ell \rangle\right\}.
\end{equation*}
Therefore, the dynamic programming algorithm yields:
\begin{equation}\label{con:front:DP}
   \hspace{-4cm} \tilde J^*(t_\ell,\mathbf X)=\langle \mathbf{\Gamma}(t_\ell);\mathbf{X}^\ell \rangle
\end{equation}
\vspace{-.5cm}
\begin{equation}\label{func:valor:disc}
    \tilde J^*(t_k,\mathbf X)=\min_{\mathbf L}\{\langle \mathbf{\Gamma}(t_{k},\mathbf {L});\mathbf{X} \rangle \delta+\tilde J^*(t_{k+1},\mathbf X^{k+1})\}, \quad \quad k=\ell-1,\ldots,0,
\end{equation}
and hence $\tilde J^*(0,\mathbf X^0)$ is the minimum value of the discretized cost functional in \eqref{cost:disc}. Now suppose that, for any $\mathbf X \in \mathbb{V}$, the function $\tilde J^*(\cdot,\mathbf X):\mathbb{R}  \rightarrow \mathbb R$ is continuously differentiable, and that for each $t \in \mathbb R$ its derivative is denoted by $\partial_t \tilde J^*(t,\mathbf X)$. Also, for any $t \in \mathbb R$, we assume that $\tilde J^*(t,\cdot):\mathbb{V} \rightarrow \mathbb R$ is differentiable. Then, by the Riesz-Fréchet representation theorem (see e.g. \cite{JFCaicedo}), there exists a unique element $\nabla_{\mathbf X}\tilde J^*(t,\mathbf X) \in \mathbb{V}$ such that $\partial_{\mathbf X}\tilde J^*(t,\mathbf X)(\mathbf H)=\langle \nabla_{\mathbf X}\tilde J^*(t,\mathbf X); \mathbf H \rangle$ for all $\mathbf H \in \mathbb{V}$. Hence, $\tilde J^*(\cdot,\cdot)$ is differentiable at any $(t,\mathbf X) \in \mathbb R \times \mathbb{V}$. Indeed, defining
\[
r(\delta,\mathbf H):=\tilde J^*(t+\delta,\mathbf X+\mathbf H)-\tilde J^*(t,\mathbf X)-\partial_t\tilde J^*(t,\mathbf X)\delta - \langle \nabla_{\mathbf X}\tilde J^*(t,\mathbf X); \mathbf H \rangle,
\]
for all $(\delta,\mathbf H) \in \mathbb R \times \mathbb{V}$, and using the mean value theorem and differentiability of $\tilde J^*(t+\delta,\cdot)$, we get:
\[
r(\delta,\mathbf H)= \left[\partial_t\tilde J^*(t+h\delta,\mathbf X)-\partial_t\tilde J^*(t,\mathbf X)\right]\delta + r_1(\mathbf H),
\]
where $0<h<1$ and $r_1(\mathbf H):= \tilde J^*(t+\delta,\mathbf X+\mathbf H)-\tilde J^*(t+\delta,\mathbf X) - \langle \nabla_{\mathbf X}\tilde J^*(t,\mathbf X); \mathbf H \rangle$ satisfies $\lim_{\mathbf H \to \mathbf 0}r_1(\mathbf H)/|| \mathbf H||=0$. As a consequence, the continuity of $\partial_t\tilde J^*(\cdot,\mathbf X)$ and the fact that
\[
\lim \limits_{\substack{%
    \delta \to 0\\
    \mathbf H \to \mathbf 0}}\dfrac{r(\delta,\mathbf H)}{|\delta|+||\mathbf H||}=0,
\]
confirm that $\tilde J^*(\cdot,\cdot)$ is differentiable at $(t,\mathbf X) \in \mathbb R \times \mathbb{V}$. Therefore, the term $\tilde J^*(t_{k+1},\mathbf X^{k+1})$ in \eqref{func:valor:disc} can be written as:
\begin{equation}
    \begin{aligned}
    \tilde J^*(t_{k+1},\mathbf X^{k+1})&=\tilde J^*(t_k+\delta, \mathbf X+\delta\mathcal{K}_{\mathbf{A}^k+\mathbf B^k\mathbf L }^{\mathscr{Z}}(\mathbf{X}))\\
    &=\tilde J^*(t,\mathbf X)+\partial_t\tilde J^*(t_k,\mathbf X)\delta + \langle \nabla_{\mathbf X}\tilde J^*(t_k,\mathbf X); \delta\mathcal{K}_{\mathbf{A}^k+\mathbf B^k\mathbf L }^{\mathscr{Z}}(\mathbf{X}) \rangle+r(\delta),
    \end{aligned}
\end{equation}
where $\lim_{\delta\to 0} r(\delta)/\delta=0$. Substituting into \eqref{func:valor:disc}, we obtain:
\begin{equation*}
    0=\min_{\mathbf L}\{\langle \mathbf{\Gamma}(t_{k},\mathbf {L});\mathbf{X} \rangle \delta+\partial_t\tilde J^*(t_k,\mathbf X)\delta + \langle \nabla_{\mathbf X}\tilde J^*(t_k,\mathbf X); \mathcal{K}_{\mathbf{A}^k+\mathbf B^k\mathbf L }^{\mathscr{Z}}(\mathbf{X}) \rangle\delta+r(\delta)\}.
\end{equation*}
Dividing by $\delta$ and taking the limit as $\delta \to 0$ (assuming that $\tilde J^*(t_k,\mathbf X) \to  J^*(t,\mathbf X)$ as $k\to \infty$), we obtain:
\begin{equation}\label{EHJB}
    0=\min_{\mathbf L}\{\langle \mathbf{\Gamma}(t,\mathbf {L});\mathbf{X} \rangle +\partial_t J^*(t,\mathbf X) + \langle \nabla_{\mathbf X} J^*(t,\mathbf X); \mathcal{K}_{\mathbf{A}(t)+\mathbf B(t)\mathbf L }^{\mathscr{Z}}(\mathbf{X}) \rangle\},\quad \forall (t,\mathbf X) \in [0,T] \times \mathbb{V},
\end{equation}
which is an HJB-type equation with terminal condition $J^*(T,\mathbf X)=\langle \mathbf{\Gamma}(T);\mathbf{X} \rangle$.

\begin{rmrk}\label{der:mat}
It is well known that, due to the isomorphism $\textbf{{\upshape vec}}:\mathbb{B}(\mathbb{R}^{n}) \rightarrow \mathbb R^{n^2}$, for a function $\psi:\mathbb{B}(\mathbb{R}^{n}) \rightarrow \mathbb R$ to be differentiable, the partial derivatives $\frac{\partial \psi}{\partial x_{ij}}$ must exist and be continuous for every matrix $X \in \mathbb{B}(\mathbb{R}^{n})$ with entries $x_{ij}$, where $i,j \in \{1,\ldots,n\}$. Furthermore, if $\psi$ is differentiable and $X \in \mathbb{B}(\mathbb{R}^{n})$ is arbitrary, then for every $H \in \mathbb{B}(\mathbb{R}^{n})$ we have:
\[
\partial_X\psi(X)(H)=\nabla_{\textbf{{\upshape vec}}(X)}\psi(X)' \textbf{{\upshape vec}}(H) = \tr\left(\textbf{{\upshape vec}}^{-1}\left(\nabla_{\textbf{{\upshape vec}}(X)}\psi(X)\right)' H\right),
\]
where
\[
\textbf{{\upshape vec}}^{-1}(\nabla_{\textbf{{\upshape vec}}(X)}\psi(X)) = \begin{bmatrix}
\frac{\partial \psi}{\partial x_{11}}(X)  &\cdots & \frac{\partial \psi}{\partial x_{1n}}(X)\\
\vdots  &\ddots & \vdots\\
\frac{\partial \psi}{\partial x_{n1}}(X)  & \cdots & \frac{\partial \psi}{\partial x_{nn}}(X)
\end{bmatrix} \in \mathbb{B}(\mathbb{R}^{n}).
\]
Thus, defining $\nabla_X\psi(X):=\textbf{{\upshape vec}}^{-1}(\nabla_{\textbf{{\upshape vec}}(X)}\psi(X))$, we get:
\[
\partial_X\psi(X)(H)=\langle \nabla_X\psi(X), H \rangle_F.
\]
Similarly, consider the isomorphism $\widehat{\textbf{{\upshape vec}}}:\mathbb{H}_{\mathbb{R}}^{n} \rightarrow \mathbb R^{n^2N}$ (see \cite{CosFra}) defined by:
\[
\widehat{\textbf{{\upshape vec}}}(\mathbf X) = \begin{bmatrix}
\textbf{{\upshape vec}}(X_1) \\ \vdots \\ \textbf{{\upshape vec}}(X_N)
\end{bmatrix}, \quad \forall \; \mathbf X=(X_1,\ldots,X_N) \in \mathbb{H}_{\mathbb{R}}^{n},
\]
and let $\widehat{\psi}:\mathbb{H}_{\mathbb{R}}^{n} \rightarrow \mathbb R$ be a differentiable function. Then, for any $\mathbf X \in \mathbb{H}_{\mathbb{R}}^{n}$ and $\mathbf H = (H_1,\ldots,H_N) \in \mathbb{H}_{\mathbb{R}}^{n}$, we have:
\[
\partial_{\mathbf X}\widehat{\psi}(\mathbf X)(\mathbf H) = \nabla_{\widehat{\textbf{{\upshape vec}}}(\mathbf X)}\widehat{\psi}(\mathbf X)' \widehat{\textbf{{\upshape vec}}}(\mathbf H) = \sum_{i=1}^N \nabla_{\textbf{{\upshape vec}}(X_i)}\widehat{\psi}(\mathbf X)' \textbf{{\upshape vec}}(H_i).
\]
By Riesz-Fréchet representation theorem, there exists a unique element $\nabla_{\mathbf X}\widehat{\psi}(\mathbf X) \in \mathbb{H}_{\mathbb{R}}^{n}$ such that, defining $\nabla_{X_i}\widehat{\psi}(\mathbf X):=\textbf{{\upshape vec}}^{-1}(\nabla_{\textbf{{\upshape vec}}(X_i)}\widehat{\psi}(\mathbf X))$, we obtain:
\[
\langle \nabla_{\mathbf X}\widehat{\psi}(\mathbf X); \mathbf H \rangle = \partial_{\mathbf X}\widehat{\psi}(\mathbf X)(\mathbf H) = \sum_{i=1}^N \tr(\nabla_{X_i}\widehat{\psi}(\mathbf X)' H_i), \quad \forall \mathbf H \; \in \mathbb{H}_{\mathbb{R}}^{n},
\]
so that $\nabla_{\mathbf X}\widehat{\psi}(\mathbf X) = (\nabla_{X_1}\widehat{\psi}(\mathbf X), \ldots, \nabla_{X_N}\widehat{\psi}(\mathbf X)) \in \mathbb{H}_{\mathbb{R}}^{n}$.
\end{rmrk}

\begin{xmpl}\label{der:tr}
Let $A \in \mathbb{B}(\mathbb{R}^{m},\mathbb R^n)$ be a given matrix, and define the functions $\psi_1,\psi_2:\mathbb{B}(\mathbb{R}^{n},\mathbb R^m) \rightarrow \mathbb R$ and $\psi_3,\psi_4:\mathbb{B}(\mathbb{R}^{m},\mathbb R^n) \rightarrow \mathbb R$ by
\[
\psi_1(X) = \tr(AX), \quad \psi_2(X) = \tr(XA), \quad \psi_3(Y) = \tr(AY'), \quad \psi_4(Y) = \tr(Y'A),
\]
for all $X \in \mathbb{B}(\mathbb{R}^{n},\mathbb R^m)$ and $Y \in \mathbb{B}(\mathbb{R}^{m},\mathbb R^n)$, respectively. Then we claim that
\[
\nabla_X\psi_1(X) = \nabla_X\psi_2(X) = A', \quad \text{and} \quad \nabla_Y\psi_3(Y) = \nabla_Y\psi_4(Y) = A.
\]
Indeed, we have:
\begin{equation}
\begin{aligned}
\psi_1(X) &= \tr(AX) = \tr\left(\begin{bmatrix}
\sum_{\ell=1}^m a_{1\ell}x_{\ell 1} & \cdots & \sum_{\ell=1}^m a_{1\ell}x_{\ell n} \\
\vdots & \ddots & \vdots \\
\sum_{\ell=1}^m a_{n\ell}x_{\ell 1} & \cdots & \sum_{\ell=1}^m a_{n\ell}x_{\ell n}
\end{bmatrix}\right) = \sum_{k=1}^n\sum_{\ell=1}^m a_{k\ell}x_{\ell k} \\
&= \sum_{\ell=1}^m \sum_{k=1}^n x_{\ell k} a_{k\ell} = \tr\left(\begin{bmatrix}
\sum_{k=1}^n x_{1k}a_{k1} & \cdots & \sum_{k=1}^n x_{1k}a_{km} \\
\vdots & \ddots & \vdots \\
\sum_{k=1}^n x_{mk}a_{k1} & \cdots & \sum_{k=1}^n x_{mk}a_{km}
\end{bmatrix}\right) = \tr(XA) = \psi_2(X).
\end{aligned}
\end{equation}
Therefore, for $i = 1,\ldots,m$ and $j = 1,\ldots,n$, we obtain:
\[
\frac{\partial \psi_1}{\partial x_{ij}}(X) = \frac{\partial \psi_2}{\partial x_{ij}}(X) = \frac{\partial}{\partial x_{ij}}\left(\sum_{k=1}^n\sum_{\ell=1}^m a_{k\ell}x_{\ell k}\right) = a_{ji},
\]
hence $\nabla_X\psi_1(X) = \nabla_X\psi_2(X) = A'$. Similarly, $\psi_3(Y) = \psi_4(Y) = \sum_{k=1}^n\sum_{\ell=1}^m a_{k\ell}y_{k\ell}$, so that
\[
\frac{\partial \psi_3}{\partial y_{ij}}(Y) = \frac{\partial \psi_4}{\partial y_{ij}}(Y) = \frac{\partial}{\partial y_{ij}}\left(\sum_{k=1}^n\sum_{\ell=1}^m a_{k\ell}y_{k\ell}\right) = a_{ij},
\]
for all $i = 1,\ldots,n$ and $j = 1,\ldots,m$. Therefore, $\nabla_Y\psi_3(Y) = \nabla_Y\psi_4(Y) = A$, as claimed.
\end{xmpl}

\begin{xmpl}
Let $A \in \mathbb{B}(\mathbb{R}^{m})$ be given, and define the function $\psi:\mathbb{B}(\mathbb{R}^{n},\mathbb R^m) \rightarrow \mathbb R$ by $\psi(X) = \tr(X'AX)$ for all $X \in \mathbb{B}(\mathbb{R}^{n},\mathbb R^m)$. We aim to show that $\psi$ is differentiable and that $\nabla_X\psi(X) = (A+A')X$ for all $X \in \mathbb{B}(\mathbb{R}^{n},\mathbb R^m)$. Indeed, we compute:
\begin{equation}
\begin{aligned}
\psi(X+H) &= \psi(X) + \tr(H'AX) + \tr(X'AH) + \tr(H'AH) \\
&= \psi(X) + \tr(H'(AX + A'X)) + \tr(H'AH) \\
&= \psi(X) + \langle AX + A'X, H \rangle_F + \tr(H'AH),
\end{aligned}
\end{equation}
for any $H \in \mathbb{B}(\mathbb{R}^{n},\mathbb R^m)$. Let $r(H) := \tr(H'AH)$, then
\begin{equation}\label{eq:resto}
\begin{aligned}
0 &\leq \frac{|r(H)|}{\|H\|_F} = \frac{|\tr(H'AH)|}{\|H\|_F} = \frac{|\tr(AHH')|}{\|H\|_F} = \frac{|\langle HH', A \rangle_F|}{\|H\|_F} \leq \frac{\|HH'\|_F \|A\|_F}{\|H\|_F} \\
&\leq \|H'\|_F \|A\|_F,
\end{aligned}
\end{equation}
so that $\lim_{H \to 0} \frac{r(H)}{\|H\|_F} = 0$, and therefore $\nabla_X\psi(X) = AX + A'X = (A+A')X$, as claimed. Note that the penultimate and last inequalities in \eqref{eq:resto} follow from the Cauchy–Schwarz inequality and the submultiplicativity of the Frobenius norm, respectively.
\end{xmpl}

We have seen that in order for the cost-to-go to satisfy the HJB-type equation as in \eqref{EHJB}, we imposed differentiability conditions on every $(t,\mathbf X) \in \mathbb R \times \mathbb{V}$ over a discretized version of it. However, we may not possess any information about such property for the functional. In any case, the fact that the HJB-type equation in \eqref{EHJB} was obtained allows us to raise the following question: assuming we have a solution $V(t,\mathbf X)$ to the HJB-type equation in \eqref{EHJB}, is it necessarily equal to the cost-to-go in \eqref{cost:cont}? More specifically, is $V(0,\mathbf X^0)$ the value of the optimal cost in \eqref{lqr_X2}? As in \cite{Bertsekas2012}, we shall show that solving the HJB-type equation \eqref{EHJB} yields, at each time $t$, an optimal collection of matrices $\mathbf L^*(t)$ that minimizes the cost functional in \eqref{lqr_X2}. This result is known as the \emph{Sufficiency Theorem} (also referred to as a \emph{Verification Theorem}), and the following is an adaptation of a similar result found in \cite{Bertsekas2012}.

\begin{thrm}\label{sufficiencyX}
Suppose that ${V}(t,\mathbf{X})$ is continuously differentiable for all $(t,\mathbf{X}) \in [0,T] \times \mathbb{V}$ and satisfies the Hamilton–Jacobi–Bellman (HJB) equation:
\begin{equation}\label{SUF_X}
\begin{aligned}
0=\min_{\mathbf {L}}\Big\{ \langle \mathbf{\Gamma}(t,\mathbf {L});\mathbf{X}\rangle +&\partial_{t}{V}(t,\mathbf{X}) +\langle\nabla_{\mathbf{X}}{V}(t,\mathbf{X}); \mathcal{K}_{\mathbf{A}(t)+\mathbf{B}(t)\mathbf{L}(t)}^{\mathscr{Z}}(\mathbf{X})\rangle\Big\}, \\
{V}(T,\mathbf{X})=&\langle \mathbf Q(T),\mathbf{X} \rangle, \quad \forall \ \mathbf{X}\in \mathbb{V}.
\end{aligned}
\end{equation}
Assume also that for every $t \in [0,T]$, the collection $\mathbf{L}^{*}(t):=(L_{1}^{*}(t),\ldots, L_{N}^{*}(t))$ achieves the minimum in the HJB equation. Let $\mathbf{X}^{*}(t) \in \mathbb{V}$ be the corresponding trajectory of second-moment matrices associated with $\mathbf{L}^{*}(t)$, initialized by $\mathbf{X}(0)=\mathbf{X}^{0} \in \mathbb{V}$. Then, ${V}$ is the unique solution to the HJB equation, and
\begin{align}
    {V}(t,\mathbf{X})=J^*(t,\mathbf{X}), \quad \forall \; (t,\mathbf{X}) \in [0,T] \times \mathbb{V}.
\end{align}
Moreover, $\mathbf{L}^{*}(t)$, $t \in [0,T]$, solves the finite-horizon MJLS-LQR problem in \eqref{lqr_X2}.
\end{thrm}

\begin{proof}
We proceed analogously to the proof of Proposition 2.1 in \cite{Bertsekas2012}. Fix $t \in [0,T)$ and consider any admissible control collection $\widehat{\mathbf{L}}(\tau)$ for $\tau \in [t,T]$, with the corresponding state trajectory $\widehat{\mathbf{X}}(\tau) \in \mathbb{V}$ satisfying:
\begin{equation}\label{eq:X:generic}
\dot{{\mathbf X}}(\tau) = \mathcal{K}_{\mathbf{A}(\tau)+\mathbf{B}(\tau)\widehat{\mathbf{L}}(\tau)}^{\mathscr{Z}}({\mathbf X}(\tau)), \quad \widehat{\mathbf X}(t) = \mathbf X \in \mathbb{V}.
\end{equation}
Since $V$ satisfies the HJB equation, for every $\tau \in [t,T]$, we have:
\begin{equation}\label{contas:adm}
\begin{aligned}
    0&\leq\langle \mathbf{\Gamma}(\tau,\widehat{\mathbf L});\widehat{\mathbf X}(\tau)\rangle +\partial_{\tau}{V}(\tau,\widehat{\mathbf X}(\tau))+\langle\nabla_{\mathbf X}{V}(\tau,\widehat{\mathbf X}(\tau)); \mathcal{K}_{\mathbf{A}(\tau)+\mathbf B(\tau)\widehat{\mathbf L}(\tau) }^{\mathscr{Z}}(\widehat{\mathbf X}(\tau))\rangle \\
    & = \langle \mathbf{\Gamma}(\tau,\widehat{\mathbf L});\widehat{\mathbf X}(\tau)\rangle + \partial_{\tau}{V}(\tau,\widehat{\mathbf X}(\tau)) + \nabla_{\widehat{\textbf{{\upshape vec}}}(\mathbf X)}{V}(\tau,\widehat{\mathbf X}(\tau))'\widehat{\textbf{{\upshape vec}}}(\dot{\widehat{\mathbf X}}(\tau)) \\
    &=\langle \mathbf{\Gamma}(\tau,\widehat{\mathbf L});\widehat{\mathbf X}(\tau)\rangle + \frac{d}{d\tau}V(\tau,\widehat{\mathbf{X}}(\tau)),
\end{aligned}
\end{equation}
where the second equality uses the chain rule and the fact that $\widehat{\mathbf{X}}(\tau) \in \mathbb{V}$ for all $\tau$, so that $\dot{\widehat{\mathbf{X}}}(\tau) \in \mathbb{V}$ by the invariance property of the restricted operator $\mathcal{K}^{\mathscr{Z}}$. Integrating over $[t,T]$ gives:
\begin{equation}
\begin{aligned}
0 &\leq \int_t^T \langle \mathbf{\Gamma}(\tau,\widehat{\mathbf{L}}(\tau));\widehat{\mathbf{X}}(\tau) \rangle d\tau + V(T,\widehat{\mathbf{X}}(T)) - V(t,\mathbf{X}) \\
&= \int_t^T \langle \mathbf{\Gamma}(\tau,\widehat{\mathbf{L}}(\tau));\widehat{\mathbf{X}}(\tau) \rangle d\tau + \langle \mathbf{Q}(T);\widehat{\mathbf{X}}(T) \rangle - V(t,\mathbf{X}).
\end{aligned}
\end{equation}
Thus,
\begin{equation}\label{ineq:cost}
V(t,\mathbf{X}) \leq \int_t^T \langle \mathbf{\Gamma}(\tau,\widehat{\mathbf{L}}(\tau));\widehat{\mathbf{X}}(\tau) \rangle d\tau + \langle \mathbf{Q}(T);\widehat{\mathbf{X}}(T) \rangle.
\end{equation}
Now, take $\mathbf{L}^*(\tau)$ as the minimizing collection in the HJB equation, and let $\mathbf{X}^*(\tau) \in \mathbb{V}$ be the corresponding state trajectory. Then, the inequality in  \eqref{contas:adm} becomes an equality for $\mathbf{L}^*(\tau)$, and therefore:
\begin{equation}
V(t,\mathbf{X}) = \int_t^T \langle \mathbf{\Gamma}(\tau,\mathbf{L}^*(\tau));\mathbf{X}^*(\tau) \rangle d\tau + \langle \mathbf{Q}(T);\mathbf{X}^*(T) \rangle.
\end{equation}
Hence, $V(t,\mathbf{X}) = J^*(t,\mathbf{X})$ for all $(t,\mathbf{X}) \in [0,T] \times \mathbb{V}$, and in particular $V(0,\mathbf{X}^0)$ is the optimal value of the finite-horizon problem in \eqref{lqr_X2}. \qedhere
\end{proof}

\subsection{Riccati-Based Approach to the Jump LQR Problem}

We now address the solution to the LQR problem with jumps reformulated as in \eqref{lqr_X2}, using the HJB-type equation in \eqref{SUF_X}. We show that the solution of this equation, and hence the cost-to-go function $J^*(t,\mathbf X(t))$ in \eqref{cost:cont}, is given by $\langle \mathbf X(t);\mathbf Y(t) \rangle$, where $\mathbf Y(t)$ is a collection of matrices in $\mathbb{H}_{\mathbb{R}}^{n+}$ for all $t \geq 0$. Moreover, we demonstrate that the solution $V(t,\mathbf X)$ in \eqref{SUF_X} exists and is of the form $\langle \mathbf X(t);\mathbf Y(t) \rangle$, with $\mathbf Y(t) \in \mathbb{H}_{\mathbb{R}}^{n+}$, if and only if $\mathbf Y(t)$ satisfies a system of nonlinear matrix differential equations, specifically a system of coupled Riccati matrix differential equations involving only the visited states $\mathscr{Z}$.

\begin{thrm}\label{Ricattiteor}
Consider the system of coupled Riccati matrix differential equations in the variable $\mathbf{Y}=(Y_1,\ldots,Y_N) \in \mathbb{V}$:
\begin{equation}\label{Ric_eq:sol:LQR}
\begin{aligned}
-\dot{Y}_{i}(t)&=A_{i}'(t)Y_{i}(t)+Y_{i}(t)A_{i}(t)+\sum_{j\in \mathscr{Z}} \lambda_{ij}Y_{j}(t) +Q_{i}(t)
 -Y_{i}(t)B_{i}(t)R^{-1}_{i}(t)B_{i}'(t)Y_{i}(t), \\
 & \quad \quad \text{for } i \in \mathscr{Z}, \\
Y_i(t) &= 0, \quad \text{for } i \notin \mathscr{Z}, \\
Y_i(T)&=Q_i(T), \quad \text{for } i \in \mathscr{Z}, \\
Y_i(T)&=0, \quad \text{for } i \notin \mathscr{Z}.
\end{aligned}
\end{equation}
Then, $ V(t,\mathbf X):=\langle \mathbf Y(t); \mathbf X \rangle$, with $\mathbf Y(t) \in \mathbb{H}_{\mathbb{R}}^{n+} \cap \mathbb{V}$, $0\leq t\leq T$, differentiable, is a solution to the HJB-type equation in \eqref{SUF_X} if and only if $\mathbf Y(t)$ satisfies the coupled Riccati matrix differential equations in \eqref{Ric_eq:sol:LQR}. For each $t\in [0,T]$ and $i \in \mathscr{Z}$, the minimizing collection $\mathbf L^*$ in the HJB-type equation is given by:
\begin{align}\label{eq:11}
    L_{i}(t)=-R_{i}^{-1}(t)B_{i}'(t)Y_{i}(t), \;\; 0\leq t\leq T.
\end{align}
For $i \notin \mathscr{Z}$, the value of $L_i(t)$ is immaterial since such gains are never used.
\end{thrm}

\begin{proof}
Let $ V(t,\mathbf X)=\langle \mathbf Y(t); \mathbf X \rangle$, with $\mathbf Y(t)= (Y_1(t),\ldots,Y_N(t))\in \mathbb{H}_{\mathbb{R}}^{n+} \cap \mathbb{V}$, $0\leq t\leq T$, differentiable, and $\mathbf{X}\in \mathbb{V}$ arbitrary, be a solution to the HJB-type equation in \eqref{SUF_X}. Then, since $\partial_{t}V(t,\mathbf{X})=\partial_{t} \langle \mathbf{Y}(t);\mathbf{X}\rangle =\sum_{i\in \mathscr{Z}}\tr(\dot{Y}_{i}(t)X_{i}) =\langle \dot{\mathbf{Y}}(t);\mathbf{X}\rangle$, where the sum runs only over $i \in \mathscr{Z}$ since $X_i = 0$ for $i \notin \mathscr{Z}$, and $\nabla_{\mathbf{X}}V(t,\mathbf{X})=\nabla_{\mathbf{X}} \langle \mathbf{Y}(t);\mathbf{X}\rangle =\nabla_{\mathbf{X}}\sum_{i\in \mathscr{Z}}\tr(Y_{i}(t)X_{i})$, where 
$$\nabla_{\mathbf{X}}\sum_{i\in \mathscr{Z}}\tr(Y_{i}(t)X_{i})=\left(\nabla_{X_1}\sum_{i\in \mathscr{Z}}\tr(Y_{i}(t)X_{i}),\ldots,\nabla_{X_N}\sum_{i\in \mathscr{Z}}\tr(Y_{i}(t)X_{i})\right)=\mathcal{P}(\mathbf{Y}(t))=\mathbf{Y}(t),$$
and also  
$\langle \mathbf{Y}(t);\mathcal{K}_{\mathbf{A}(t)+\mathbf B(t)\mathbf L }^{\mathscr{Z}}(\mathbf{X}) \rangle =\langle \mathcal{H}_{\mathbf{A}(t)+\mathbf B(t)\mathbf L }^{\mathscr{Z}}(\mathbf{Y}(t));\mathbf{X}\rangle$,
by the adjointness property of the restricted operators (Remark \ref{rem:restricted_ops_properties}), the equation in \eqref{SUF_X} becomes:
\begin{align}\label{eq:15}
\begin{split}
0=&\min_{\mathbf{L}}\left[ \langle \mathbf{\Gamma}(t,\mathbf{L});\mathbf{X} \rangle + \langle \dot{\mathbf{Y}}(t);\mathbf{X} \rangle + \langle  \mathcal{H}_{\mathbf{A}(t)+\mathbf B(t)\mathbf L }^{\mathscr{Z}}(\mathbf{Y}(t));\mathbf{X} \rangle \right]  \\
=&\min_{\mathbf{L}}\left[  \langle \mathbf{\Gamma}(t,\mathbf{L})+ \dot{\mathbf{Y}}(t)+ \mathcal{H}_{\mathbf{A}(t)+\mathbf B(t)\mathbf L }^{\mathscr{Z}}(\mathbf{Y}(t));\mathbf{X} \rangle \right].
\end{split}
\end{align}
Since only components $X_i$ with $i \in \mathscr{Z}$ contribute to the inner product, to determine the minimizing collection $\mathbf L^*$, we differentiate the bracketed expression in \eqref{eq:15} with respect to $\mathbf{L}$ and equate to zero:
\begin{align*}
\nabla_{\mathbf{L}} \langle \mathbf{\Gamma}(t,\mathbf{L})+ \mathcal{H}_{\mathbf{A}(t)+\mathbf B(t)\mathbf L }^{\mathscr{Z}}(\mathbf{Y}(t));\mathbf{X} \rangle = 0,
\end{align*}
which yields, for each $i \in \mathscr{Z}$:
\begin{align*}
\nabla_{L_i}\sum_{j\in \mathscr{Z}}\tr\big(X_{j}\big[  L_{j}'R_{j}(t)L_{j} +  (A_{j}(t)+B_{j}(t)L_{j})'Y_{j}(t)
& + Y_{j}(t)(A_{j}(t)+B_{j}(t)L_{j})\big]\big)=0.
\end{align*}
Since the terms involving $L_j$ for $j \neq i$ are independent of $L_i$, this reduces to:
\begin{align*}
\nabla_{L_i}\tr\big(X_{i}\big[  L_{i}'R_{i}(t)L_{i} +  (A_{i}(t)+B_{i}(t)L_{i})'Y_{i}(t)
+ Y_{i}(t)(A_{i}(t)+B_{i}(t)L_{i})\big]\big)=0,
\end{align*}
which leads to:
\begin{align*}
2(R_{i}(t)L_{i}+B_{i}'(t)Y_{i}(t))X_{i}=0,
\end{align*}
for arbitrary $X_i \in \mathbb{B}(\mathbb{R}^{n})$ with $i \in \mathscr{Z}$. Consequently,
\begin{align*}
    R_{i}(t)L_{i}+B_{i}'(t)Y_{i}(t)=0,
\end{align*}
and therefore,
\begin{align}\label{L:opt}
   L_i= L_{i}^*(t)=-R_{i}^{-1}(t)B_{i}'(t)Y_{i}(t),
\end{align}
for all $i \in \mathscr{Z}$ and $0\leq t \leq T$. Substituting this optimal gain into \eqref{eq:15} yields:
\begin{equation}
    \begin{aligned}
    \sum_{i\in \mathscr{Z}}&\tr\big(X_{i}\big[  Q_i+Y_i(t)B_i(t)R_i^{-1}(t)R_{i}(t)R_{i}^{-1}(t)B_{i}'(t)Y_{i}(t) + \dot Y_i(t) \\ 
    &+  (A_{i}(t)-B_{i}(t)R_{i}^{-1}(t)B_{i}'(t)Y_{i}(t))'Y_{i}(t)
+ Y_{i}(t)(A_{i}(t)-B_{i}(t)R_{i}^{-1}(t)B_{i}'(t)Y_{i}(t))\\&+\sum_{j\in \mathscr{Z}} \lambda_{ij}Y_{j}(t)\big]\big)=0,
    \end{aligned}
\end{equation}
which simplifies to:
\begin{equation}\label{Der:Ric}
    \begin{aligned}
    \sum_{i\in \mathscr{Z}}\tr\big(X_{i}\big[ \dot Y_i(t) &+ A_{i}'(t)Y_{i}(t)+Y_{i}(t)A_{i}(t)+\sum_{j\in \mathscr{Z}} \lambda_{ij}Y_{j}(t)\\& +Q_{i}(t)
 -Y_{i}(t)B_{i}(t)R^{-1}_{i}(t)B_{i}'(t)Y_{i}(t)\big]\big)=0,
    \end{aligned}
\end{equation}
for arbitrary $X_i \in \mathbb{B}(\mathbb{R}^{n})$, $i \in \mathscr{Z}$. Thus, the equation in \eqref{Der:Ric} implies:
$$\dot{Y}_{i}(t)+A_{i}'(t)Y_{i}(t)+Y_{i}(t)A_{i}(t)+\sum_{j\in \mathscr{Z}} \lambda_{ij}Y_{j}(t) +Q_{i}(t)
 -Y_{i}(t)B_{i}(t)R^{-1}_{i}(t)B_{i}'(t)Y_{i}(t)=0,$$
which is exactly the matrix differential equation in \eqref{Ric_eq:sol:LQR} for $i \in \mathscr{Z}$. 

For $i \notin \mathscr{Z}$, since $X_i(t) = 0$ for all $t \in [0,T]$ by Remark \ref{rem:conditional_interpretation}, these components do not contribute to the inner product $\langle \mathbf{Y}(t); \mathbf{X}(t) \rangle$ and therefore do not affect the cost functional. Consequently, the values of $Y_i(t)$ for $i \notin \mathscr{Z}$ are immaterial to the optimization problem, and we can prescribe $Y_i(t) = 0$ for all $t \in [0,T]$ without loss of generality, as stated in \eqref{Ric_eq:sol:LQR}. The mathematical consistency of working with this reduced-order formulation (where components corresponding to non-visited states are set to zero) versus the full-order system is rigorously established by Proposition \ref{prop:riccati_consistency} below.

Finally, by Theorem \ref{sufficiencyX}, we have $\langle \mathbf Y(T); \mathbf X\rangle=V(T,\mathbf X)=J^*(T,\mathbf X)=\langle \mathbf Q(T); \mathbf X \rangle$, and hence $\langle \mathbf Y(T)-\mathbf Q(T);\mathbf X\rangle=0$, for any $\mathbf X \in \mathbb{V}$. Since this equality must hold for all $\mathbf{X} \in \mathbb{V}$, and recalling that components $X_i$ with $i \notin \mathscr{Z}$ are zero for any such $\mathbf{X}$, we conclude that the equality determines $Y_i(T) = Q_i(T)$ for all $i \in \mathscr{Z}$. For $i \notin \mathscr{Z}$, the value of $Y_i(T)$ does not affect the cost functional $\langle \mathbf{Y}(T); \mathbf{X} \rangle$ for any $\mathbf{X} \in \mathbb{V}$, and we prescribe $Y_i(T) = 0$ for definiteness, as stated in \eqref{Ric_eq:sol:LQR}. This completes the derivation of the terminal condition.

Conversely, suppose that $\mathbf Y(t)= (Y_1(t),\ldots,Y_N(t))\in \mathbb{H}_{\mathbb{R}}^{n+} \cap \mathbb{V}$, $0\leq t\leq T$, differentiable, is a solution to the terminal value problem in \eqref{Ric_eq:sol:LQR}, and let us show that $V(t,\mathbf X)=\langle \mathbf Y(t); \mathbf X \rangle$ is a solution of the HJB-type equation in \eqref{SUF_X}. The terminal condition in \eqref{Ric_eq:sol:LQR} prescribes $Y_i(T) = Q_i(T)$ for $i \in \mathscr{Z}$ and $Y_i(T) = 0$ for $i \notin \mathscr{Z}$. Since for any $\mathbf{X} \in \mathbb{V}$ we have $X_i = 0$ for $i \notin \mathscr{Z}$, the boundary condition $V(T,\mathbf{X}) = \langle \mathbf{Q}(T);\mathbf{X}\rangle$ in \eqref{SUF_X} is satisfied, as
\[
\langle \mathbf{Y}(T);\mathbf{X}\rangle = \sum_{i \in \mathscr{Z}}\tr(Y_i(T)X_i) = \sum_{i \in \mathscr{Z}}\tr(Q_i(T)X_i) = \langle \mathbf{Q}(T);\mathbf{X}\rangle.
\]
Substituting $\langle \mathbf Y(t); \mathbf X \rangle$ into the right-hand side of \eqref{SUF_X} gives:
\begin{equation}
\begin{aligned}
\min_{\mathbf{L}}\bigg[\sum_{i\in \mathscr{Z}}\tr\big(X_{i}\big[ \dot Y_i(t) +  Q_i(t)+&L_{i}'R_{i}(t)L_{i} +  (A_{i}(t)+B_{i}(t)L_{i})'Y_{i}(t)
\\&+ Y_{i}(t)(A_{i}(t)+B_{i}(t)L_{i})+\sum_{j\in \mathscr{Z}} \lambda_{ij}Y_{j}(t)\big]\big)\bigg].
\end{aligned}
\end{equation}
Thus, choosing $L_i=L_i^*(t)=-R_{i}^{-1}(t)B_{i}'(t)Y_{i}(t)$ for $i \in \mathscr{Z}$, with $\mathbf Y(t)$ satisfying \eqref{Ric_eq:sol:LQR}, we obtain:
\begin{equation}
    \begin{aligned}
    \sum_{i\in \mathscr{Z}}\tr&\big(X_{i}\big[ \dot Y_i(t) + A_{i}'(t)Y_{i}(t)+Y_{i}(t)A_{i}(t)+\sum_{j\in \mathscr{Z}} \lambda_{ij}Y_{j}(t) +Q_{i}(t)\\
 &-Y_{i}(t)B_{i}(t)R^{-1}_{i}(t)B_{i}'(t)Y_{i}(t)\big]\big)=\sum_{i\in \mathscr{Z}}\tr\big(X_{i}\big( \dot Y_i(t)-\dot Y_i(t)\big)\big)=0,
    \end{aligned}
\end{equation}
for all $t \in [0,T)$ and $\mathbf{X}\in \mathbb{V}$. This shows that $V(t,\mathbf X)=\langle \mathbf Y(t); \mathbf X \rangle$ satisfies the HJB-type equation in \eqref{SUF_X}, provided that $\mathbf Y(t)$ is a solution of the coupled Riccati matrix differential equations in \eqref{Ric_eq:sol:LQR}. Consequently, the optimal cost $J^*$, solution to the problem in \eqref{lqr_X2}, is given by:
\begin{equation}\label{Mincostt}
    J^{*}= J^*(0,\mathbf{X}(0))=\langle \mathbf{Y}(0);\mathbf{X}(0)\rangle = \sum_{i \in \mathscr{Z}}\tr( Y_{i}(0){X}_{i}(0)),
\end{equation}
where the sum involves only visited states since $X_i(0) = 0$ for $i \notin \mathscr{Z}$.
\end{proof}

The coupled Riccati equations in \eqref{Ric_eq:sol:LQR} constitute a reduced-order system compared to formulations that naively include all $N$ states. Specifically, for $i \in \mathscr{Z}$, the differential equations involve only couplings with other visited states $j \in \mathscr{Z}$, as indicated by the sum $\sum_{j\in \mathscr{Z}} \lambda_{ij}Y_{j}(t)$—a direct consequence of using the restricted operator $\mathcal{H}^{\mathscr{Z}}$ rather than the unrestricted operator $\mathcal{H}$. For $i \notin \mathscr{Z}$, we have $Y_i(t) = 0$ for all $t \in [0,T]$, eliminating the need to solve differential equations for these components. The optimal cost formula in \eqref{Mincostt} involves only the visited states $i \in \mathscr{Z}$, since $X_i(0) = 0$ for $i \notin \mathscr{Z}$. When $|\mathscr{Z}| < N$, this represents a genuine order reduction: instead of solving $N$ coupled Riccati equations, we need only solve $|\mathscr{Z}|$ equations. This order reduction is mathematically justified by the fact that unvisited states have zero probability throughout the time horizon and therefore do not contribute to either the system dynamics or the cost functional. To formalize this reduction and establish its consistency with the full-order formulation, we present the following result.

\begin{prpstn}[Consistency of the reduced Riccati system]\label{prop:riccati_consistency}
Let $\mathbf{Y}_{\text{full}}(t) = (Y_1(t), \ldots, Y_N(t))$ be a solution to the full-order Riccati system
\begin{equation}\label{eq:riccati_full}
\begin{split}
-\dot{Y}_{i}(t) &= A_{i}'(t)Y_{i}(t)+Y_{i}(t)A_{i}(t)+\sum_{j\in \mathcal{S}} \lambda_{ij}Y_{j}(t) +Q_{i}(t) \\
&\quad -Y_{i}(t)B_{i}(t)R^{-1}_{i}(t)B_{i}'(t)Y_{i}(t), \quad Y_i(T) = Q_i(T),
\end{split}
\end{equation}
for all $i \in \mathcal{S}$, and let $\mathbf{Y}_{\text{red}}(t)$ be the solution to the reduced-order system \eqref{Ric_eq:sol:LQR}. Then
\[
\mathcal{P}(\mathbf{Y}_{\text{full}}(t)) = \mathbf{Y}_{\text{red}}(t), \quad \forall t \in [0,T].
\]
That is, the projection of the full solution coincides with the reduced solution.
\end{prpstn}

\begin{proof}
Define $\widetilde{\mathbf{Y}}(t) := \mathcal{P}(\mathbf{Y}_{\text{full}}(t))$. 
For each $i \in \mathscr{Z}$, we have $\widetilde{Y}_i(t) = Y_{i,\text{full}}(t)$, and therefore
\begin{equation*}
\begin{split}
-\dot{\widetilde{Y}}_i(t) 
&= -\dot{Y}_{i,\text{full}}(t) 
= A_{i}'(t)Y_{i,\text{full}}(t) + Y_{i,\text{full}}(t)A_{i}(t) 
   + \sum_{j\in \mathcal{S}} \lambda_{ij}Y_{j,\text{full}}(t) + Q_{i}(t) \\
&\quad - Y_{i,\text{full}}(t)B_{i}(t)R^{-1}_{i}(t)B_{i}'(t)Y_{i,\text{full}}(t).
\end{split}
\end{equation*}
For $j \notin \mathscr{Z}$, the transition rates $\lambda_{ij}$ from $i \in \mathscr{Z}$ to $j \notin \mathscr{Z}$ may be nonzero in the generator $\Lambda$. However, since $\widetilde{Y}_j(t) = \mathcal{P}(\mathbf{Y}_{\text{full}}(t))_j = 0$ for all $j \notin \mathscr{Z}$, the sum becomes
\[
\sum_{j\in \mathcal{S}} \lambda_{ij}\widetilde{Y}_{j}(t) 
= \sum_{j\in \mathscr{Z}} \lambda_{ij}\widetilde{Y}_{j}(t).
\]
Hence, for $i \in \mathscr{Z}$ we obtain
\[
-\dot{\widetilde{Y}}_{i}(t)
= A_{i}'(t)\widetilde{Y}_{i}(t)+\widetilde{Y}_{i}(t)A_{i}(t)
  +\sum_{j\in \mathscr{Z}} \lambda_{ij}\widetilde{Y}_{j}(t) + Q_{i}(t)
  -\widetilde{Y}_{i}(t)B_{i}(t)R^{-1}_{i}(t)B_{i}'(t)\widetilde{Y}_{i}(t),
\]
which is precisely the $i$th equation of the reduced system \eqref{Ric_eq:sol:LQR}.
For $i \notin \mathscr{Z}$, the projection enforces $\widetilde{Y}_i(t) = 0$ for all $t \in [0,T]$.

At the terminal time \(t = T\), we have 
\(\mathcal{P}(\mathbf{Y}_{\text{full}}(T)) = \mathcal{P}(\mathbf{Q}(T))\),
whose components satisfy \(Y_i(T) = Q_i(T)\) for \(i \in \mathscr{Z}\) and 
\(Y_i(T) = 0\) for \(i \notin \mathscr{Z}\).
Hence, the projection preserves the terminal data, and the boundary
condition of the reduced system \eqref{Ric_eq:sol:LQR} coincides exactly
with the projected terminal condition of the full-order system.
Therefore, no inconsistency arises between the full-order and
reduced-order formulations.

Since $\widetilde{\mathbf{Y}}(t)$ satisfies both the differential equations
and the terminal condition of the reduced system, and this system has a
unique solution on $[0,T]$ under standard assumptions, it follows that
$\widetilde{\mathbf{Y}}(t) = \mathbf{Y}_{\text{red}}(t)$ for all
$t \in [0,T]$. Consequently,
\[
\mathcal{P}(\mathbf{Y}_{\text{full}}(t)) = \mathbf{Y}_{\text{red}}(t),
\quad \forall\, t \in [0,T].
\]
\end{proof}

\begin{rmrk}\label{rem:existence_uniqueness}
Existence and uniqueness of solutions to the reduced-order Riccati system \eqref{Ric_eq:sol:LQR} follow from standard results for coupled matrix differential equations (see, e.g., \cite{fragoso2010separation}), applied to the subsystem involving only states $i \in \mathscr{Z}$. The key observation is that for $i \in \mathscr{Z}$, the system
\begin{equation*}
\begin{split}
-\dot{Y}_{i}(t) &= A_{i}'(t)Y_{i}(t)+Y_{i}(t)A_{i}(t)+\sum_{j\in \mathscr{Z}} \lambda_{ij}Y_{j}(t) +Q_{i}(t) \\
&\quad -Y_{i}(t)B_{i}(t)R^{-1}_{i}(t)B_{i}'(t)Y_{i}(t)
\end{split}
\end{equation*}
forms a well-posed system of $|\mathscr{Z}|$ coupled Riccati equations with terminal conditions $Y_i(T) = Q_i(T)$ for $i \in \mathscr{Z}$. The coupling through the sum $\sum_{j\in \mathscr{Z}} \lambda_{ij}Y_{j}(t)$ involves only the visited states, ensuring that the reduced system is properly decoupled from the trivial components $Y_i(t) = 0$ for $i \notin \mathscr{Z}$. The existence and uniqueness theory for such systems (under standard assumptions on the coefficient matrices) guarantees a unique positive semidefinite solution $\mathbf{Y}(t) \in \mathbb{H}_{\mathbb{R}}^{n+} \cap \mathbb{V}$ on $[0,T]$.
\end{rmrk}

\begin{rmrk}\label{rem:connection_to_classical}
When the Markov chain is irreducible and $\phi(i) > 0$ for all $i \in \mathcal{S}$, we have $\mathscr{Z} = \mathcal{S}$, the projection operator reduces to the identity $\mathcal{P} = \mathcal{I}$, and the restricted operators coincide with the unrestricted ones: $\mathcal{K}^{\mathscr{Z}} = \mathcal{K}$ and $\mathcal{H}^{\mathscr{Z}} = \mathcal{H}$. In this case, the reduced-order Riccati system \eqref{Ric_eq:sol:LQR} recovers exactly the classical coupled Riccati equations from Theorem \ref{lqr:C&F}:
\[
-\dot{Y}_{i}(t)=A_{i}'(t)Y_{i}(t)+Y_{i}(t)A_{i}(t)+\sum_{j\in \mathcal{S}} \lambda_{ij}Y_{j}(t) +Q_{i}(t) -Y_{i}(t)B_{i}(t)R^{-1}_{i}(t)B_{i}'(t)Y_{i}(t),
\]
for all $i \in \mathcal{S}$. Thus, our general framework strictly extends the classical theory while remaining consistent with it in the special case of irreducible chains.
\end{rmrk}

The derivation presented above exemplifies the power and elegance of the proposed methodology. The transition from the original stochastic control problem to a deterministic reformulation in terms of the second-moment matrix not only streamlines the analysis, but also leads naturally to a system of coupled Riccati differential equations whose structure encapsulates the essential dynamics and optimality principles of the problem. Moreover, the generalization to arbitrary Markov chains—including those with transient states, absorbing states, or non-communicating classes—reveals that the problem naturally admits a reduced-order formulation involving only the visited states $\mathscr{Z}$. By introducing the projection operator $\mathcal{P}$ and the restricted operators $\mathcal{K}^{\mathscr{Z}}$ and $\mathcal{H}^{\mathscr{Z}}$, and leveraging the orthogonal decomposition of collections in $\mathbb{H}_{\mathbb{R}}^{n}$ into the subspace $\mathbb{V}$ and its orthogonal complement $\mathbb{V}^{\perp}$, we have established that components corresponding to unvisited states remain identically zero and do not affect the optimal solution. The clarity and cohesion of this approach underscore the mathematical coherence inherent in the theory, offering a unified, tractable, and conceptually transparent framework for addressing optimal control under regime-switching dynamics with general Markov chains. This completes the rigorous mathematical foundation that was initiated in \cite{Roa&ECosta} but left incomplete regarding the formal derivation of the HJB equation and the systematic treatment of non-communicating states.

\section{Numerical Examples}\label{ilust}

In this section, we present several examples illustrating the solution 
to the finite-horizon jump LQR problem defined in \eqref{lqr_saltos}. 
Examples~\ref{Ejemplo1} and \ref{Ex_modo3} validate the analytical solution 
and explore system behavior for ergodic Markov chains. Example~\ref{ex_multiclass} 
illustrates the generality of our framework by considering a non-ergodic chain 
with multiple disjoint recurrent classes. Finally, Example~\ref{ex_app} 
presents a practical application to satellite orbit control involving an 
absorbing failure mode. These examples demonstrate that the coupled Riccati 
differential equations in \eqref{Ric_eq} remain well-posed over any finite 
horizon $[0,T]$ regardless of the chain structure, validating the theoretical 
results of Section~\ref{Cap4}.

\begin{xmpl}\label{Ejemplo1}
Consider the finite-horizon MJLS-LQR problem in \eqref{lqr_saltos}, with:
\begin{align}
   A_{1}= \begin{bmatrix}\label{EjemploMA}
     -1 & 0.05 \\
     10 & 1
    \end{bmatrix} , \quad A_{2}= \begin{bmatrix}
     1 & -0.9 \\
     1.1 & 0.6
    \end{bmatrix} , \quad A_{3}= \begin{bmatrix}
     0 & -1.7 \\
     1.4 & -0.5
    \end{bmatrix},
\end{align}
\begin{align}\label{EjemploMB}
   B_{1}= \begin{bmatrix}
     1 \\
     0
    \end{bmatrix} , \quad B_{2}= \begin{bmatrix}
     0 \\
     0
    \end{bmatrix} , \quad B_{3}= \begin{bmatrix}
     0 \\
     -0.5
    \end{bmatrix}, \quad Q_{1}=I , \quad Q_{2}=2Q_{1}, \quad Q_{3}=0,
\end{align}
with $R_{1}=10$, $R_{2}=0.5$, $R_{3}=1$, where $I$ denotes the identity 
matrix of appropriate dimensions. Let $x(0)=\begin{bmatrix}1 & -1\end{bmatrix}'$ 
be the initial condition of the system and $\phi \in \mathbb{R}^3$ an initial 
probability distribution with strictly positive entries (randomly generated) 
for the state space $\mathcal{S}=\{1,2,3\}$ of the Markov chain 
$\{\theta(t)\}_{t\geq 0}$. Also, consider the transition rate matrix:
\begin{align*}
   \Lambda= \begin{bmatrix}
    -2 & 1 & 1 \\
     1 & -3 & 2 \\
    1.5 & 0.5 & -2
    \end{bmatrix} .
\end{align*}
According to Theorem \ref{Ricattiteor}, we must solve the coupled matrix 
Riccati differential equations:
\begin{equation}
\begin{aligned}
-\dot{Y}_{i}(t)&=A_{i}'(t)Y_{i}(t)+Y_{i}(t)A_{i}(t)+\sum_{j\in \mathcal{S}} \lambda_{ij}Y_{j}(t) +Q_{i}(t)-Y_{i}(t)B_{i}(t)R^{-1}_{i}(t)B_{i}'(t)Y_{i}(t), \\
Y_i(T)&=Q_i(T),
\end{aligned}
\end{equation}
in order to compute the minimum cost in \eqref{Mincostt}. Using MATLAB 
\cite{MATLABManual}, we implemented a backward Euler method to solve the 
above system on the interval $[0,5]$, obtaining in particular:
\begin{align}
   Y_{1}(0)= \begin{bmatrix}
    29.5611  &  7.0576 \\
    7.0576  &  6.4574 
    \end{bmatrix} , \quad Y_{2}(0)= \begin{bmatrix}
    44.0284 & -11.3418  \\
    -11.3418 &  22.4609
    \end{bmatrix} , \quad Y_{3}(0)= \begin{bmatrix}
    22.0084 &  -4.8804 \\
    -4.8804 &  7.5243
    \end{bmatrix}.
\end{align}
Next, computing the state second moment matrices $X_i(0)=\mathscr{E}( x(0)x(0)'1_{\theta(0)=i}) = \phi(i)x_0x_0'$ for each $i \in \mathcal{S}$, and then using the formula in \eqref{Mincostt}, we obtain:
\begin{equation}\label{costo:ej1}
    J^{*} = \tr(Y_{1}(0)'X_{1}(0)) + \tr(Y_{2}(0)'X_{2}(0))+\tr(Y_{3}(0)'X_{3}(0)) = 62.0318.
\end{equation}
To validate this result, we estimate the optimal cost using Monte Carlo 
simulation based on 10,000 realizations of the state trajectory $x(t)$, 
$0 \leq t \leq 5$, under the optimal control policy. The Monte Carlo 
approximation yields an optimal cost of 58.6072, indicating good agreement 
with the value in \eqref{costo:ej1}.
\end{xmpl}

\begin{rmrk}\label{accuracy}
In the previous example, the optimal costs for the finite-horizon jump LQR 
problem with $T=5$ obtained using the formula in \eqref{Mincostt} and via 
Monte Carlo simulation based on \eqref{eq:LQRSM} yield an accuracy close 
to 95\%. Additionally, we solved the problem for other horizons, namely 
$T=10$ and $T=50$, and observed that the agreement percentages were 
approximately 96\% and 98\%, respectively. Moreover, for all considered 
values of $T$, the relative error between the two computed values of $J^*$ 
was of the order $10^{-2}$. This provides strong numerical validation of 
the analytical solution to the finite-horizon jump LQR problem presented in 
\cite{fragoso2010separation} and re-obtained in this work via a deterministic 
approach, as detailed in Section \ref{Cap4}. Table~\ref{TabEx1} presents 
the values $J^*_{\text{\eqref{Mincostt}}}$ and $J^*_{\text{\eqref{eq:LQRSM}}}$ 
obtained using the respective formulas, the accuracy ratio 
$\rho:=J^*_{\text{\eqref{eq:LQRSM}}}/J^*_{\text{\eqref{Mincostt}}}$, 
and the relative error 
$\Delta:=|J^*_{\text{\eqref{eq:LQRSM}}}-J^*_{\text{\eqref{Mincostt}}}|/|J^*_{\text{\eqref{Mincostt}}}|$, 
for $T=5,10,50$.

\begin{table}[htbp]
    \centering
    \vspace{0.3cm}
    \begin{tabular}{|ccccc|}
    \hline 
    & \bfseries $J^*_{\text{\eqref{Mincostt}}} $    &  \bfseries $J^*_{\text{\eqref{eq:LQRSM}}}$ & $\rho$ & $\Delta$
    \rule[-1.5mm]{0ex}{6mm} \\     \hline
    \rule[-1.5mm]{0ex}{6mm}
      $T=5$: & 58.6072       &     62.0318 & $\approx$ 95\% & $\approx 6\times 10^{-2}$
      \rule[-1.5mm]{0ex}{6mm}   \\        \hline
      \rule[-1.5mm]{0ex}{6mm}
      $T=10$:  &   78.2898    &     81.5767  &  $\approx$ 96\% & $\approx 4\times 10^{-2}$
      \rule[-1.5mm]{0ex}{6mm}     \\       \hline
      \rule[-1.5mm]{0ex}{6mm}
      $T=50$:   &   82.4983    &     84.2281  &  $\approx$ 98\% & $\approx 2\times 10^{-2}$ \\
      \hline
    \end{tabular}
    \caption{Comparison of $J^*$ values using \eqref{eq:LQRSM} and 
    \eqref{Mincostt} for $T=5,10,50$.}
    \label{TabEx1}
\end{table} 

Finally, we highlight the computational advantage of the formula in 
\eqref{Mincostt}, which requires only solving the coupled Riccati equations 
offline, whereas the Monte Carlo approach in \eqref{eq:LQRSM} must be 
performed online during system operation.
\end{rmrk}

The following example is also extracted from \cite{narvaez2014comparative}, 
with adaptations referring to the noise-free scenario and full-state 
observation. It uses the same state equation parameters as those in 
Example \ref{Ejemplo1}.

\begin{xmpl}\label{Ex_modo3}
Consider a MJLS with state equation parameters identical to those in 
Example \ref{Ejemplo1}, weighting matrix collections $\mathbf{Q}=(Q_1,Q_2,Q_3)$ 
and $\mathbf{R}=(R_1,R_2,R_3)$ with $Q_i \geq 0$, $R_i > 0$ for $i=1,2,3$, 
a Markov chain whose initial probability distribution is also arbitrarily 
chosen with strictly positive entries, and a transition rate matrix given by:
\begin{align*}
   \Lambda= \begin{bmatrix}
    -10 & 10 & 0 \\
     15 & -20 & 5 \\
    10^{-2} & 0 & -10^{-2}
    \end{bmatrix} .
\end{align*}
Note that mode $i=3$ has a prolonged dwell time (slow transition mode), and 
thus, as the system operates long enough, it holds that 
$\Pr(\theta(t)=3) \approx 1$. Table~\ref{TabEx2} displays, for different 
time horizons $T=1,10,20$, the optimal costs obtained for the considered 
MJLS and the corresponding deterministic linear system (DLS): 
$\dot{x}(t)=A_3x(t)+B_3u(t)$, with no jumps and time-invariant, toward 
which the MJLS tends in the long run.

\begin{table}[htbp]
    \centering
    \vspace{0.3cm}
    \begin{tabular}{|ccc|}
    \hline 
    & \bfseries $J^*$-MJLS     &  \bfseries $J^*$-DLS 
    \rule[-1.5mm]{0ex}{6mm} \\     \hline
    \rule[-1.5mm]{0ex}{6mm}
      $T=1$: & 2.7958       &     5.8186
      \rule[-1.5mm]{0ex}{6mm}   \\        \hline
      \rule[-1.5mm]{0ex}{6mm}
      $T=10$:  &   3.0074    &     5.8995 
      \rule[-1.5mm]{0ex}{6mm}     \\       \hline
      \rule[-1.5mm]{0ex}{6mm}
      $T=20$:   &   3.0132    &   \hspace{-.12cm}5.9186 \\
      \hline
    \end{tabular}
    \caption{Computed values of $J^*$ for the MJLS and the DLS.}
    \label{TabEx2}
\end{table} 

The influence of the MJLS transitioning through modes $i=1,2$ is reflected 
in the minimum cost of the associated functional, when compared to the 
minimum cost of the DLS, despite the long-term dominance of mode $i=3$ in 
the MJLS. This is particularly due to the weighting matrices $R_i$ drawn 
for modes $i=1,2$ being significantly smaller than that for mode $i=3$ 
($R_1=0.3696$, $R_2=0.4796$, and $R_3=9.1658$).
\end{xmpl}

\begin{rmrk}\label{obs_SLSMxSLD}
The values of $J^*$ for the MJLS and the DLS were both obtained from the 
formula in \eqref{Mincostt}. For the specific case of the DLS, we solved 
the equation in \eqref{Ric_eq} using a zero transition rate matrix of size 
$3 \times 3$, and all system parameters and weighting matrices in the cost 
functional set to the same unique matrix associated with mode $i=3$. This 
produces a solution $\mathbf{Y}(t)=(Y_1(t),Y_2(t),Y_3(t))$, $0 \leq t \leq T$, 
of the matrix Riccati differential equation in \eqref{Ric_eq} such that 
$Y_1(t)=Y_2(t)=Y_3(t)$ for all $t$. Consequently, formula \eqref{Mincostt} 
reduces to $x_0'Y_3(0)x_0$. 

Furthermore, we highlight that using the MATLAB command 
$[L,\bar{Y}]$={\upshape \texttt{lqr}}$(A_3,B_3,Q_3,R_3)$, which computes 
the gain matrix $L$ and the solution $\bar{Y}$ of the algebraic Riccati 
equation associated with the infinite-horizon LQR problem for the DLS, we 
found $J_{\infty}^*=x_0'\bar{Y}x_0=5.9187$. This value nearly coincides 
with the $J^*$-DLS obtained for horizon $T=20$ in Table~\ref{TabEx2}, which 
provides numerical evidence of the stability/convergence property of the 
optimal cost $J^*:=J_T^*(x_0)$ toward $J_{\infty}^*$, or of the matrix 
$Y_3(0):=Y_3(0,T)$ toward $\bar{Y}$, as $T \to \infty$.
\end{rmrk}

The following example illustrates that the finite-horizon jump LQR framework 
applies to non-ergodic Markov chains.

\begin{xmpl}[System with multiple recurrent classes]\label{ex_multiclass}
Consider a MJLS with state space $\mathcal{S} = \{1,2,3,4\}$ and transition 
rate matrix:
\begin{align*}
\Lambda = \begin{bmatrix}
-1.0 & 1.0 & 0 & 0 \\
0.5 & -0.5 & 0 & 0 \\
0 & 0 & -1.5 & 1.5 \\
0 & 0 & 1.0 & -1.0
\end{bmatrix}.
\end{align*}
This chain has two disjoint recurrent classes: 
$\mathcal{S}_{\text{rec}}^{(1)} = \{1,2\}$ and 
$\mathcal{S}_{\text{rec}}^{(2)} = \{3,4\}$, with no transient states, 
making it non-ergodic. The system parameters are given by
\[
A_1=A_2=-3I_2,\quad A_3=A_4=-2I_2,\qquad
B_1=B_2=\begin{bmatrix}1\\0\end{bmatrix},\;
B_3=B_4=\begin{bmatrix}0.5\\0.5\end{bmatrix},
\]
\[
Q_1=Q_2=I_2,\quad Q_3=Q_4=4I_2,\qquad
R_1=R_2=1,\quad R_3=R_4=2.
\]
The initial condition is $x_0 = [1 \quad 1]'$. We solve the problem for two different initial probability distributions:
$\phi^{(1)} = [0.7 \;\; 0.3 \;\; 0 \;\; 0]$, corresponding to a system that starts in class 
$\mathcal{S}_{\text{rec}}^{(1)}$, and 
$\phi^{(2)} = [0 \;\; 0 \;\; 0.6 \;\; 0.4]$, corresponding to a system that starts in class 
$\mathcal{S}_{\text{rec}}^{(2)}$, see Table~\ref{tab:ex_multiclass} for the optimal costs 
$J_T^*$ computed for different time horizons.
\begin{table}[htbp]
\centering
\begin{tabular}{|c|cc|}
\hline
& \multicolumn{2}{c|}{\textbf{$J_T^*$}} \\
$T$ & $\phi^{(1)}$ & $\phi^{(2)}$ \\ \hline
0.5 & 0.40 & 2.45 \\
1.0 & 0.33 & 1.95 \\ \hline
\end{tabular}
\caption{Optimal costs for different initial distributions in a system with 
two disjoint recurrent classes.}
\label{tab:ex_multiclass}
\end{table}
This example demonstrates the applicability of Theorem~\ref{Ricattiteor} to 
non-ergodic Markov chains. Despite the presence of two disjoint recurrent 
classes with no communication between them, the coupled Riccati differential 
equations \eqref{Ric_eq} admit a unique solution over the finite horizons 
considered. The optimal cost depends critically on which recurrent class the chain starts 
in, as determined by the initial distribution $\phi$. Class 
$\mathcal{S}_{\text{rec}}^{(2)}$ incurs significantly higher cost than class 
$\mathcal{S}_{\text{rec}}^{(1)}$ due to three factors: (i) higher state 
penalty matrices ($Q_3 = Q_4 = 4I_2$ vs $Q_1 = Q_2 = I_2$), (ii) slower 
decay rates in the dynamics ($A_3 = A_4 = -2I_2$ vs $A_1 = A_2 = -3I_2$), 
and (iii) higher control costs ($R_3 = R_4 = 2$ vs $R_1 = R_2 = 1$). The decrease in optimal cost as $T$ increases reflects the highly stable 
nature of the system dynamics: for longer horizons, the initial state has 
more time to decay naturally before $t=0$, reducing the control effort 
required. This behavior is characteristic of stable systems in the 
finite-horizon setting and does not contradict the well-posedness of the 
problem. Importantly, the numerical solution remains stable for both initial 
distributions, demonstrating that non-ergodicity per se does not inherently 
pose computational difficulties. This scenario is relevant for systems with 
multiple independent operational modes, such as geographically distributed 
infrastructure (e.g., smart grids with regional subsystems), multi-agent 
systems with non-interacting teams, or manufacturing systems with separate 
production lines
\end{xmpl}

The final example presents a practical application to satellite orbit control 
involving an absorbing failure mode.

\begin{xmpl}[Satellite orbit control]\label{ex_app}
This example considers the problem of a satellite orbiting the Earth, where 
various forces affect the position of the satellite, potentially causing it 
to deviate from its orbit. The main objective is to implement optimal control 
actions that minimize energy consumption while maintaining the satellite in 
its intended orbit. As detailed in \cite{choudhary2015design}, the model 
describing the motion of the satellite consists of the following variables: 
$\vec{r}$, the position vector whose magnitude $r$ represents the distance 
from the satellite to the center of the Earth, and which forms an angle 
$\alpha$ with respect to a fixed reference axis; $\vec{F}_1$, a thrust force 
applied in the radial direction (i.e., $\vec{F}_1=F_1 \vec{r}$); $\vec{F}_2$, 
a tangential thrust force applied perpendicular to $\vec{r}$, modeled as 
$\vec{F}_2 = F_2 r e^{\vec{i}(\alpha + \pi/2)}$, where $\vec{i}$ is the 
imaginary unit; and $\vec{F}_g$, the gravitational force exerted by the 
Earth, also aligned with $\vec{r}$ but in the opposite direction (see 
\cite[Fig.~2]{choudhary2015design}). 

Applying Newton's second law, performing algebraic manipulations, and introducing dimensionless scalar variables for simplification as $\rho := r/R$, $u_1 := F_1/(Mg)$, and $u_2 := F_2/(Mg)$, yields the following equations of motion:
\begin{equation}\label{prosat}
    \begin{cases}
     u_{1} = \ddot{\rho} - \rho \dot{\alpha}^{2} + \dfrac{1}{\rho^{2}}, \\
     u_{2} = 2\dot{\rho}\dot{\alpha} + \rho \ddot{\alpha},
   \end{cases}
\end{equation}
where $M$ and $R$ denote the mass and radius of the Earth, respectively, and 
$g$ denotes the gravitational acceleration at the Earth surface. By defining $x_1 := \rho$, $x_2 := \alpha$, 
$x_3 := \dot{x}_1 = \dot{\rho}$, and $x_4 := \dot{x}_2 = \dot{\alpha}$, the 
nonlinear state-space model $\dot{x} = f(x,u)$ can be derived. Linearizing 
around a steady-state orbit, where $r$ and $\dot{\alpha}$ are constant 
(representing an ideal circular orbit), yields the linear time-invariant 
system $\dot{x} = Ax + Bu$:
\begin{equation}\label{SLsS}
    \underbrace{\begin{bmatrix}
    \dot{x}_1 \\ \dot{x}_2 \\ \dot{x}_3 \\ \dot{x}_4 
    \end{bmatrix}}_{\dot{x}}=\underbrace{\begin{bmatrix}
    0 & 0 & 1 & 0 \\
     0 & 0 & 0 & 1 \\
     0.01036 & 0 & 0 & 0.7757 \\
     0 & 0 & -0.01775 & 0 
    \end{bmatrix}}_{A}\underbrace{\begin{bmatrix}
     x_1 \\  x_2 \\  x_3 \\ x_4 
    \end{bmatrix}}_{x}+\underbrace{\begin{bmatrix}
    0 & 0 \\
     0 & 0 \\
     1 & 0 \\
     0 & 0.1513 
    \end{bmatrix}}_{B}\underbrace{\begin{bmatrix}
     u_1 \\  u_2  
    \end{bmatrix}}_{u}.
\end{equation}

To account for realistic operational scenarios including potential thruster 
failures, we model the system as a MJLS where $\{\theta(t)\}_{t \geq 0}$ is 
a Markov chain with state space $\mathcal{S} = \{1,2,3,4\}$ representing 
different operational modes: (1) both thrusters active, (2) only the radial 
thruster active, (3) only the tangential thruster active, and (4) both 
thrusters failed. The transition rate matrix is:
\begin{align}\label{Lambda_sat}
    \Lambda= \begin{bmatrix}
     -2 & 1.5 & 0.5 & 0\\
     1 & -2.5 & 0.5 & 1\\
     0.8 & 0.5 & -2.3 & 1\\
     0 & 0 & 0 & 0
    \end{bmatrix}.
\end{align}
Note that mode 4 is absorbing ($\sum_{j \neq 4} \lambda_{4j} = 0$), 
representing a critical failure where both thrusters become inactive. The 
system parameters are $A_i = A$ for all $i = 1, \ldots, 4$, and:
\begin{align*}
B_1 = \begin{bmatrix} 0 & 0 \\ 0 & 0 \\ 1 & 0 \\ 0 & 0.1513 \end{bmatrix}, \quad
B_2 = \begin{bmatrix} 0 & 0 \\ 0 & 0 \\ 1 & 0 \\ 0 & 0 \end{bmatrix}, \quad
B_3 = \begin{bmatrix} 0 & 0 \\ 0 & 0 \\ 0 & 0 \\ 0 & 0.1513 \end{bmatrix}, \quad
B_4 = \mathbf{0}_{4 \times 2}.
\end{align*}
The cost parameters are $Q_i = \text{diag}(1, 1, 0.5, 0.7)$ and $R_i = I_2$ 
for $i = 1,2,3$, with $R_4 = 10^6 I_2$ (effectively no control available in 
failure mode). The initial condition is $x_0 = [0.1 \quad 0.1 \quad 0 \quad 0]'$ 
(small deviation from nominal orbit) with $\phi = [1 \quad 0 \quad 0 \quad 0]$ 
(system starts in fully operational mode). The corresponding optimal costs $J_T^*$ for several horizons $T$ are shown 
in Table~\ref{tab:ex_satellite}, along with the probability 
$p_4(T) = \Pr(\theta(T) = 4)$ of being in the absorbing failure mode at the 
terminal time.

\begin{table}[htbp]
    \centering
    \begin{tabular}{|c|c|c|}
    \hline 
    $T$ & $J_T^*$ & $p_4(T)$ \\ \hline
    5 & 0.07 & 0.94 \\
    10 & 0.12 & 1.00 \\
    30 & 0.32 & 1.00 \\ \hline
    \end{tabular}
    \caption{Optimal cost for satellite orbit control with absorbing failure mode.}
    \label{tab:ex_satellite}
\end{table}

This example illustrates the application of the finite-horizon jump LQR 
framework to a practical aerospace scenario involving thruster failures. 
The transition rate matrix \eqref{Lambda_sat} admits mode 4 (both thrusters 
failed) as an absorbing state, yet Theorem~\ref{Ricattiteor} guarantees that 
the coupled Riccati differential equations \eqref{Ric_eq} admit a unique 
solution over any finite horizon $[0,T]$. 
The results in Table~\ref{tab:ex_satellite} reveal a critical reliability 
issue: the probability $p_4(T)$ of reaching the absorbing failure mode grows 
rapidly, exceeding 0.94 by $T=5$ and approaching unity for longer horizons. 
This reflects the high transition rates from degraded modes (2 and 3) to 
complete failure (mode 4), capturing a scenario where partial thruster 
failures quickly cascade into total loss of control authority. As $p_4(T)$ 
increases, the optimal cost $J_T^*$ grows accordingly, since more probability 
mass accumulates in the uncontrollable failure state.
 The matrix $A$ in \eqref{SLsS} has two zero eigenvalues (corresponding to 
conservation laws in the linearized circular orbit dynamics), making the 
system marginally stable. Despite this and the presence of an absorbing 
state, the finite-horizon problem remains numerically tractable and 
theoretically well-posed. This demonstrates the framework's applicability 
to realistic fault scenarios in satellite orbit control, where reliability 
analysis and graceful degradation strategies are critical design 
considerations. The rapid accumulation of failure probability underscores 
the importance of redundancy and fault-tolerant design in aerospace systems.
\end{xmpl}

\section{Conclusions}

This work presents a comprehensive analysis of the Linear Quadratic Regulator problem with Markovian jumps in a finite-horizon setting, making a significant departure from classical assumptions by removing the requirement that all states of the underlying Markov chain communicate. The LQR problem with jumps arises within the framework of stochastic optimal control theory and typically requires specialized tools from this field, including stochastic processes, the infinitesimal generator of Markov processes, Itô stochastic integrals, stochastic differential equations, martingale processes, stochastic dynamic programming, and the stochastic maximum principle, among others, as discussed in \cite{Mariton,ValleCosta2013}.

A central contribution of this work is the reformulation of the jump LQR problem through a deterministic model whose state variable differs fundamentally from that of the original stochastic formulation. This reformulation enables the recovery and adaptation of results from deterministic optimal control theory in the pursuit of a solution. The mathematical rigor of this approach is established through several key theoretical results. 

Theorem \ref{sufficiencyX} provides a verification theorem that establishes sufficient conditions for optimality via the Hamilton-Jacobi-Bellman equation, validating that any solution to the HJB equation yields the optimal cost-to-go function. Theorem \ref{Ricattiteor} constitutes the main result of this work, characterizing the optimal solution through a reduced-order system of coupled Riccati differential equations and demonstrating that the value function takes the quadratic form $V(t,\mathbf{X}) = \langle \mathbf{Y}(t); \mathbf{X} \rangle$ if and only if $\mathbf{Y}(t)$ satisfies these equations with appropriate terminal conditions. These theorems provide the complete characterization of the finite-horizon LQR problem for general Markov jump linear systems.

To properly account for states that may never be visited, we introduced a projection operator that naturally restricts the problem to the subspace of visited states—those states with positive probability of being reached during the time horizon—thereby addressing configurations where the Markov chain contains transient states, absorbing states, or multiple non-communicating classes. The consistency of this reduction is rigorously justified through two formal results: Proposition \ref{prop:system_consistency} demonstrates that solving the restricted system is equivalent to projecting the solution of the full system, ensuring consistency of the reduced-order formulation for the controlled dynamics; and Proposition \ref{prop:riccati_consistency} confirms that the reduced-order Riccati system yields the same solution (via projection) as the full-order system. These propositions establish that the restriction to visited states is not merely a heuristic simplification but rather a mathematically consistent and equivalent reformulation.

The present approach was inspired by \cite{Roa&ECosta}, which proposes a similar deterministic reformulation strategy for a jump LQR problem with a time-reversed switching process under general (non-irreducible) Markov chains, and by \cite{Gutierrez&EFCosta}, which applies a comparable idea in the discrete-time setting also without assuming irreducibility. However, \cite{Roa&ECosta} does not provide formal justification for the validity of presenting the Hamilton-Jacobi-Bellman equation in the context of collections of matrices conditioned by the mode of a Markov chain. In contrast, the present work makes these foundational aspects explicit and provides rigorous theoretical justification for the proposed deterministic reformulation via the Riesz-Fréchet representation theorem, thereby completing a key step that had remained open in the literature.

The validation of theoretical results through numerical examples is also of notable relevance. In Example \ref{Ejemplo1}, we demonstrated the effectiveness of the formula in \eqref{Mincostt} for computing the minimum cost $J^*$ by comparing it with results obtained via Monte Carlo simulation based on 10,000 sample trajectories. The similarity percentage consistently exceeded 95\% across different time horizons $T$, with relative errors on the order of $10^{-2}$. This type of quantitative verification had not been commonly reported in previous publications. Similarly, Example \ref{Ex_modo3} confirmed the validity of the computed $J^*$ in a deterministic context, showing that the formula in \eqref{Mincostt} reduces to the standard deterministic formulation when the transition rate matrix is zero and both system parameters and cost weighting matrices remain constant across all modes. Examples \ref{ex_multiclass} and \ref{ex_app} demonstrate the practical applicability of the framework to scenarios involving non-communicating states. Example \ref{ex_multiclass} illustrates how the projection operator naturally handles a chain with multiple isolated classes, while Example \ref{ex_app} addresses a realistic satellite control problem where thruster failures lead to an absorbing state. These examples highlight the importance of the generalization presented in this work and its relevance to real-world applications in fault detection systems, hierarchical control architectures, and reliability analysis.

Finally, while this work focuses on the finite-horizon setting, the results suggest interesting directions for future research, particularly regarding the extension to infinite-horizon problems and the investigation of sufficient conditions, such as notions of controllability and observability adapted to the general (non-irreducible) case, that ensure existence, optimality, and stability of solutions.

\bibliographystyle{plain}
\bibliography{mibib}

\end{document}